\documentclass[opre,nonblindrev]{informs3} 

\OneAndAHalfSpacedXI 

\usepackage{endnotes}
\let\footnote=\endnote

\usepackage{natbib}
 \bibpunct[, ]{(}{)}{,}{a}{}{,}%
 \def\bibfont{\small}%
 \def\newblock{\ }%

 \usepackage{bm}
 \usepackage{mathabx}
 \usepackage{subfigure}  
 \usepackage{mdwlist}
 \usepackage{cases}
 \usepackage{empheq}

\TheoremsNumberedThrough     
\ECRepeatTheorems

\EquationsNumberedThrough    

\begin{document}


\RUNAUTHOR{Pineda and Morales}

\RUNTITLE{Impact of Imbalance Costs on Stochastic Unit Investment}

\TITLE{Modeling the Impact of Imbalance Costs on Generating Expansion of Stochastic Units}

\ARTICLEAUTHORS{%
\AUTHOR{Salvador Pineda}
\AFF{University of Copenhagen, \EMAIL{s.pineda@math.ku.dk}} 
\AUTHOR{Juan M. Morales}
\AFF{Technical University of Denmark, \EMAIL{jmmgo@imm.dtu.dk}}
} 

\ABSTRACT{%
The imbalance costs incurred by a stochastic power producer due to forecast production errors have a significant impact on its total profit and therefore, such an impact needs to be taken into account when evaluating investment decisions. In this paper, we propose a modeling framework to analyze the effect of these imbalance costs on optimal generating expansion decisions of stochastic units. The proposed model is cast as a mathematical program with equilibrium constraints, which allows the explicit representation of both the day-ahead and balancing market-clearing mechanisms. We use the proposed framework to investigate the effect of two paradigmatic market designs on investment decisions: a day-ahead market that is cleared following a conventional cost merit-order principle, and an ideal market-clearing procedure that determines day-ahead dispatch decisions accounting for their impact on balancing operation costs. The variability throughout the planning horizon of the expected stochastic power production and demand level in the day-ahead market is modeled via a scenario set. Likewise, the uncertainty pertaining to their corresponding forecast errors, which are to be settled in the balancing market, is also characterized through scenarios. The main features and results of the proposed models are discussed using an illustrative two-node example and a more realistic 24-node case study.  
}%


\maketitle

%


\section{Introduction}

As renewable energy sources approach grid parity, it is widely believed that stochastic power producers should trade their electricity production just as conventional generators do. Under these conditions, a substantial share of the revenue of stochastic power producers derives from their participation in the balancing market \citep{Fabbri, matevosyan2006minimization, pinson2007trading, morales2010short}. This is so because, when it comes to trading in the day-ahead market, stochastic power producers cannot perfectly predict their future power output 
and as a result, they must turn to the balancing market to cope with their forecast errors.

The need for stochastic power producers to trade part of its production in the balancing market makes their profit intrinsically dependent on system flexibility, as only those generating units that are flexible enough to increase or decrease their power output at short notice are actually eligible to provide balancing energy. For this reason, the amount and characteristics of the flexible resources available at the real-time operation stage are likely to have a significant impact on stochastic producers' profit. Consequently, investments in new stochastic generating units are to consider not only the locations with the highest capacity factors, but also those with the most competitive access to balancing power. Likewise, the design of the electricity market may also have, by extension, an important effect on stochastic capacity investments, inasmuch as the market organization may condition the availability of flexible resources at the balancing operation of the power system \citep{weber2010adequate,borggrefe2011balancing,winkler2012market}.

In this paper, we propose a modeling framework that allows us to investigate the effect of the electricity market desing on the optimal investment decisions of a power producer that accounts for the projected revenues obtained by the new stochastic generating units in both the day-ahead and the balancing markets. The proposed model belongs to the family of mathematical programs with equilibrium constraints (MPECs), which have proved to be especially useful in the area of investment models for transmission and generation capacity expansion \citep{murphy2005generation, Garces, nanduri2009generation, wang2009strategic, wogrin2011generation, Morales}, as they naturally capture the subordination of operation outcomes to expansion decisions. More particularly, the equilibrium constraints in our investment model represent the sequential settlement of the day-ahead and the balancing markets under one of the two following market-design paradigms:
\begin{enumerate}
  \item \emph{Conventional market clearing} (ConvMC), according to which the clearing of the day-ahead market does \emph{not} anticipate the need for flexible capacity to balance the system in real time. This clearing mechanism is representative of market designs where there is no provision of reserve capacity prior to or in concurrence with the determination of the day-ahead energy dispatch.  
  \item \emph{Stochastic market clearing} (StocMC), in which the clearing of the day-ahead market \emph{does} account for the need for balancing capacity in real time by jointly optimizing the day-ahead dispatch and the subsequent balancing operation of the power system. It should be noticed that the stochastic market-clearing mechanism endogenously determines the optimal amount of reserve capacity to be procured in concurrence with the day-ahead energy dispatch \citep{bouffard2008stochastic,pritchard2010single, Morales2012}. This market-clearing mechanism is therefore representative of market designs where the provision of reserve capacity is carried out according to an ideal operating reserve demand curve.  
\end{enumerate}

The resulting optimization problems are \emph{stochastic} MPECs that consider the variability regarding both the expected stochastic power production and the expected demand level, as well as the uncertainty corresponding to the forecast errors of these two random variables. Generally, a discrete representation of variable and uncertain data, usually referred to as \emph{scenario} set, is constructed by sampling the probability distributions of the random variables \citep{Garces, wang2009strategic, wogrin2011generation, baringo2011wind}. Likewise, scenario reduction techniques are normally applied in order to trim down the computational burden of the proposed model while keeping most of the stochastic information embedded in the original probability distributions \citep{dupavcova2003scenario, morales2009scenario}.

In short, the contributions of this paper are twofold. First, we develop a modeling framework to explicitly account for the impact of imbalance costs on the optimal generation expansion decisions of a stochastic power producer. To the best of our knowledge, this feature makes the resulting investment model the first of its kind. Secondly, we analyze the effect of the market design on stochastic capacity expansion by considering two paradigmatic market-clearing mechanisms, ConvMC and StocMC, that differ in whether or not future balancing energy needs are taken into account when clearing the day-ahead electricity market.

The remainder of this paper is organized as follows. In Section \ref{SectionFormulation}, the two generation expansion models of a stochastic power producer are described and formulated as stochastic mathematical programs with equilibrium constraints. Section \ref{SectionUncertainty} presents the procedure to generate the scenario set that characterizes the uncertainty pertaining to stochastic production and demand level both in the day-ahead and the balancing stages. The main features of the proposed model are discussed in Section \ref{SectionExample} using a two-node illustrative example. Besides, a 24-node case study is presented in Section \ref{SectionCaseStudy} to analyze the impact of imbalance costs on investment decisions in a more realistic set up. Finally, Section \ref{SectionConclusions} concludes the paper.

\section{Model description and formulation} \label{SectionFormulation}

Both generation and transmission expansion planning problems have drawn interest from the scientific community for several decades now \citep{luss1982operations,latorre2003classification}. Moreover, most of the existing methodologies have also been revisited and adapted to the new competitive environment established in most power systems worldwide \citep{murphy2005generation,  wang2009strategic, Garces,    wogrin2011generation, baringo2011wind}. A number of the proposed methods to determine optimal expansion decisions are cast as bilevel programming problems due to its suitability to model the actions and reactions of a non-cooperative game between two agents with different objective functions \citep{bard1998practical}. 

\begin{figure}[htbp]
	\centering		\includegraphics{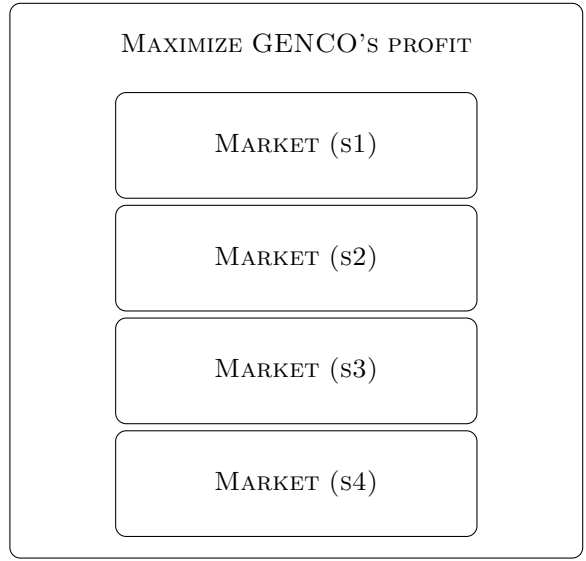}
	\caption{Bilevel generating expansion problem of a power producer participating in an electricity market}
	\label{fig:FigInvDA}
\end{figure}

More specifically, the bilevel optimization models presented in \cite{wogrin2011generation} and \cite{baringo2011wind} determine the optimal generation expansion decisions of generating companies (GenCo) participating in electricity markets. While in \cite{wogrin2011generation}, the GenCo has to decide the optimal expansion strategy of its conventional and dispatchable thermal-based generating units, \cite{baringo2011wind} propose an optimization problem that yields the optimal generation expansion decisions corresponding to stochastic generating units, such as wind farms. Both bilevel optimization models have, though, a very similar mathematical structure. The objective function of the upper level in both cases involves the maximization of the expected profit obtained from the new generating units to be built. Likewise, the lower-level problems represent the impact of the generating expansion decisions on the market outcomes for different states of the power system. Figure \ref{fig:FigInvDA} illustrates the bilevel mathematical structure of the above mentioned formulations, where \textsc{s1}, \textsc{s2}, \textsc{s3}, and \textsc{s4} represent four possible conditions of the power system characterized by, for example, demand level, stochastic production, network contingencies, etc.

The general mathematical formulation of this type of bilevel optimization models to determine investment decisions of stochastic power producers participating in electricity markets is presented below:  

\begin{subequations}\label{Bilevel1Stage}
\begin{align}
\hspace{-35mm} \underset{\bm{U}}{{\rm Max}} \quad & \sum_{s} \pi_{s} \left( \widehat{\bm{w}}^T_{s} diag(\bm{q}^T\bm{U}) \bm{\lambda}_{s} \right) - \bm{i}^T\bm{U}\bm{1}  \label{B1S_upperOF} \hspace{55mm} \\
{\rm s.t.} \quad & \bm{U}\bm{1} \leq \bm{1} \label{B1S_upperC1}
\end{align}
\vspace{-8mm}
\begin{empheq}[left={\bm{\lambda}_{s} \in \text{arg}}\empheqlbrace,right=\empheqrbrace \forall s.]{align}
\underset{\bm{p}_{s},\bm{\delta}_{s}}{{\rm  Min}} \quad & \bm{c}^T\bm{p}_{s}  \label{B1S_lowerOF} \\
 {\rm s.t.} \quad & \underline{\bm{p}_s} \leq \bm{p}_{s} \leq \widebar{\bm{p}_{s\hspace{1mm}}}: \underline{\bm{\alpha}_s}, \widebar{\bm{\alpha}_s} \label{B1S_lowerC1} \\
& -\widebar{\bm{f}_s} \leq \bm{B}\bm{A}\bm{\delta}_{s} \leq \widebar{\bm{f}_s}: \underline{\bm{\beta}_s}, \widebar{\bm{\beta}_s} \label{B1S_lowerC2} \\
& \bm{p}_s +   diag(\bm{q}^T\bm{U})\widehat{\bm{w}}_{s} + \bm{A}^T\bm{B}\bm{A}\bm{\delta}_{s} = \bm{0}: \bm{\lambda}_s \label{B1S_lowerC3} 
\end{empheq}
\end{subequations}

In this formulation and throughout the rest of this paper, capital bold letters represent matrices, while lower case bold letters symbolize vectors. Likewise, $\bm{1}$ corresponds to an all-ones vector of appropriate dimension, and $diag(\cdot)$ provides a square matrix with the elements of a given vector on the main diagonal. Dual variables are indicated at the corresponding equations following a colon. Additional notation for this formulation is given in Table \ref{tab:Notation1}.

\begin{table}[htbp]
	\centering \caption{Notation of bilevel generating expansion problem \eqref{B1S_upperOF}--\eqref{B1S_lowerC3} }
		\begin{tabular}{c l}		
		\hline		
		\multicolumn{2}{l}{Numbers and sets}\\		
		\hspace{5mm} $n_B$ & Number of nodes in the system\\	
		\hspace{5mm} $n_L$ & Number of lines in the system\\
		\hspace{5mm} $n_W$ & Number of generating expansion projects\\					\hspace{5mm} $s$   & Scenario index for market conditions\\ 
		\multicolumn{2}{l}{Parameters}\\
		\hspace{5mm}$\bm{A}$ & $n_L\times n_B$ line-to-node incidence matrix \\	
		\hspace{5mm}$\bm{B}$ & $n_L\times n_L$ diagonal matrix of line susceptances (p.u.) \\	
		\hspace{5mm}$\bm{c}$ & $n_B\times1$ vector of bid and offer prices  (\$/MWh) \\	
		\hspace{5mm}$\widebar{\bm{f}}_{s}$ & $n_L\times1$ vector of line capacities in scenario $s$ (MW) \\
		\hspace{5mm}$\bm{i}$ & $n_W\times1$ vector of annualized investment cost (\$)\\
		\hspace{5mm}$\underline{\bm{p}_{s}}$ & $n_B\times1$ vector of minimum dispatch limits in scenario $s$ (MW) \\
		\hspace{5mm}$\widebar{\bm{p}_{s\hspace{1mm}}}$ & $n_B\times1$ vector of maximum dispatch limits in scenario $s$ (MW) \\
		\hspace{5mm}$\bm{q}$ & $n_W\times1$ vector of project capacities (MW)\\
		\hspace{5mm}$\widehat{\bm{w}}_{s}$ & $n_B\times1$ vector of stochastic power production in scenario $s$ (p.u.) \\		
		\hspace{5mm}$\pi_s$ & Probability of occurrence of scenario $s$ \\
		\multicolumn{2}{l}{Variables}\\
		\hspace{5mm}$\bm{U}$ & $n_W\times n_B$ matrix of binary variables of generating expansion decisions \\	
		\hspace{5mm}$\bm{p}_{s}$ & $n_B\times1$ vector of dispatched quantities in scenario $s$ (MWh)\\
		\hspace{5mm}$\bm{\delta}_{s}$ & $n_B\times1$ vector of voltage angles in scenario $s$ (rad)\\
		\hspace{5mm}$\bm{\lambda}_{s}$ & $n_B\times1$ vector of nodal prices in scenario $s$ (\$/MWh)\\	
		\hline	
\end{tabular}\label{tab:Notation1}\end{table}

Without loss of generality, it is assumed in this formulation that generation and consumption of electricity do not occur at the same node, taking $\bm{p}_{s}$ positive values for generation nodes, and negative values for consumption nodes. Furthermore, investment decisions $\bm{U}$ are assumed to be binary, since typical constructions of new units are given in discrete and relatively large amounts to take advantage of economies of scale, being $\bm{U}(w,b)$ equal to 1 if stochastic generating unit $w$ of capacity $\bm{q}(w)$ is located at node $b$, and 0 otherwise. 

The upper-level problem \eqref{B1S_upperOF}--\eqref{B1S_upperC1} maximizes the expected profit of the stochastic power producer over all plausible scenarios $s$ that may occur during the planning horizon minus the investment costs of the undertaken projects $\bm{i}^T\bm{U}\bm{1}$. The profit for each scenario $s$ is computed as the per-unit power production in that particular scenario $\widehat{\bm{w}}^T_{s}$ times the capacity installed at each node  $diag(\bm{q}^T\bm{U})$ times the locational marginal prices $\bm{\lambda}_{s}$ provided by the lower-level market-clearing problem. Note that this objective function only contains products of binary variables and continuous variables, which can be easily linearized \citep{floudas1995nonlinear}. Additionally, constraint \eqref{B1S_upperOF} in the upper level problem enforces that each project can be placed at only one node in the system. Observe that while the stochastic producer can behave in a strategic manner by locating new generating units at those locations that will lead to the highest revenue, its selling offer price is equal to its marginal cost, which is normally very close or equal to zero. 

The electricity price at which the wind power production is paid is computed in each of the lower-level problems modeling an energy-only market-clearing procedure for each scenario $s$. Objective function \eqref{B1S_lowerOF} minimizes the production cost computed as $\bm{c}^T\bm{p}_{s}$. Minimum power outputs and capacities of generating units, as well as fix demand levels are imposed through constraint \eqref{B1S_lowerC1}. Constraint \eqref{B1S_lowerC2} enforces the transmission capacity limits. Finally, the balance equation at each node of the system is imposed through \eqref{B1S_lowerC3}, where power flows are computed according to a DC representation of the transmission network. Note that the injection of new stochastic generating units is determined as the capacity installed at each node $diag(\bm{q}^T\bm{U})$ times the per-unit power production for each particular scenario $\widehat{\bm{w}}^T_{s}$. Finally, the electricity price corresponds to the dual variable $\bm{\lambda}_s$ of the balance equation \eqref{B1S_lowerC3}. Observe that a perfectly competitive market, where producers offer to sell their whole technically available capacity at marginal cost, is considered here. Moreover, the demand is considered inelastic and can only be shed at a very high cost. 

Alternatively to the use of Karush-Kuhn-Tucker conditions, bilevel optimization problem \eqref{B1S_upperOF}-\eqref{B1S_lowerC3} can be solved by replacing the linear lower-level problems by their corresponding primal and dual constraints plus the strong duality theorem \citep{motto2005mixed}. Note that the primal-dual formulation of the lower-level problems avoids the use of an usually large number of additional binary variables required to linearize the complementarity equations corresponding to the KKT conditions. 

While investment model \eqref{Bilevel1Stage} successfully captures the impact of capacity investment decisions on market outcomes, it fails to recognize that electricity markets nowadays comprise, at least, two different trading floors of especial relevance for a stochastic power producer, namely, the day-ahead and the balancing markets. On the one hand, the day-ahead market is cleared between 24 and 36 hours in advance, thus facing a significant level of uncertainty pertaining to the generation and demand levels. On the other hand, the balancing market occurs very close to operation and deals with forecast errors by re-dispatching flexible generating units. Even though the total trading volume for balancing energy is usually very small compared to the total energy volume exchanged in the day-ahead market, the ``forced'' participation of a stochastic power producer in the balancing market may have a significant impact on its incomes. In fact, results in \cite{Fabbri, matevosyan2006minimization, pinson2007trading} show that these imbalance costs can represent between 10\% and 20\% of the total profit achieved by stochastic power producers. Consequently, the trading in the balancing market should be accounted for to determine investment decisions in new stochastic generating units. With this aim in mind, we expand problem \eqref{Bilevel1Stage} by adding a new level or stage to represent the sequential clearing process of the day-ahead and the balancing markets.

In this paper, we analyze the impact of two paradigmatic market-clearing model designs on the investment decisions of stochastic power producers. In the first market model, the day-ahead and the balancing markets are independently cleared. That is, the conditional expected power production of stochastic units is included as a zero-cost production offer in the day-ahead market clearing, and dispatched quantities are determined according to a cost-merit order criterion that takes network constraints into consideration. Subsequently, the offers for up and down regulation provided by flexible generating units are cleared at the balancing stage in order to counteract forecast errors corresponding to stochastic units and demand levels in the most economical manner possible. This model will be referred to hereinafter as \emph{Conventional Market Clearing (ConvMC)}. On the other hand, the so-called \emph{Stochastic Market Clearing (StochMC)} provides day-ahead quantities and prices accounting for the plausible future balancing needs in the system according to the stochastic characterization of the uncertain parameters and the flexibility of the system \citep{pritchard2010single, Morales2012}. In doing so, StochMC may dispatch some generating units out of merit to obtain a day-ahead generation schedule that minimizes the expected total system operation costs, i.e., including both the day-ahead and the balancing stages. Note also that, for this model, the day-ahead dispatch of the stochastic units may be different from the forecast value. 

\begin{figure}[htbp]
\centering
\subfigure[Conventional market clearing (ConvMC)]{\includegraphics {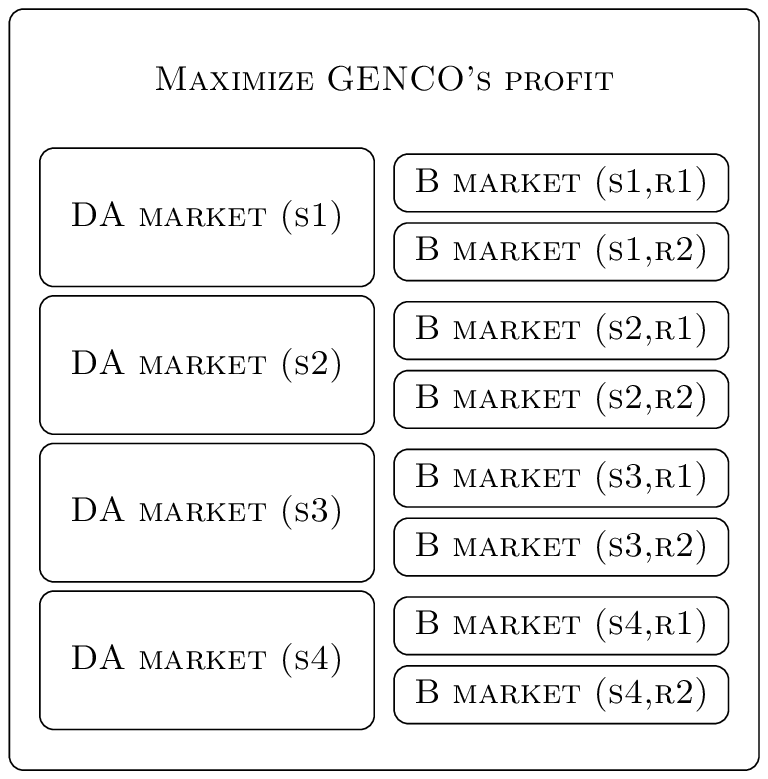}\label{fig:FigInvConv}}
\subfigure[Stochastic market clearing (StochMC)]{\includegraphics {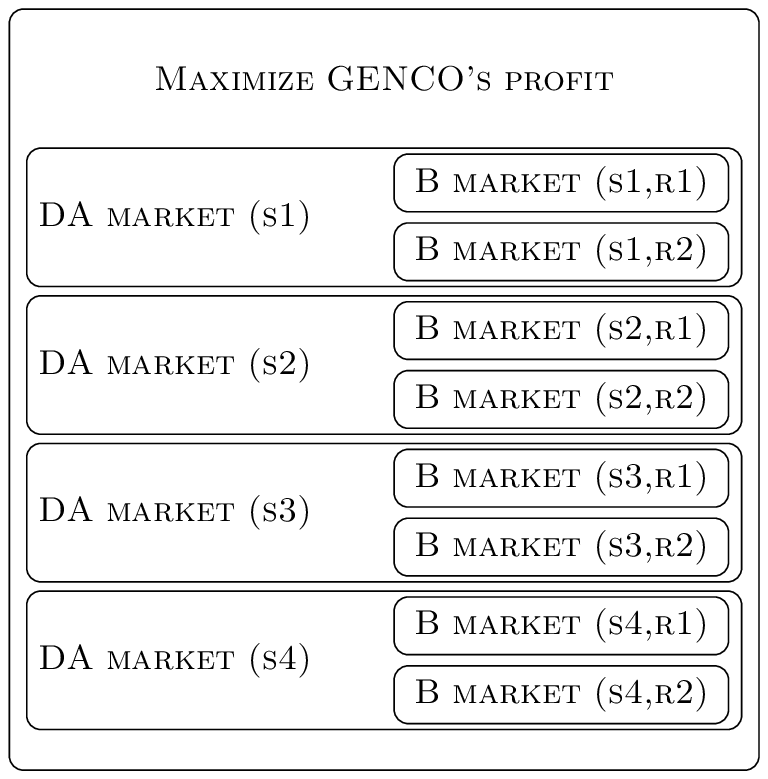}\label{fig:FigInvStoc}}
\caption{Generating expansion problem of stochatic producers considering both day-ahead and balancing markets}\label{fig:FigInv2markets}
\end{figure}

Figure \ref{fig:FigInv2markets} illustrates the two proposed mathematical frameworks to analyze the impact of imbalance costs on generating expansion decisions of stochastic power producers. For the sake of illustration, only two different scenarios at the balancing stage ({\scshape r1,r2}) are considered for each of the four plausible outcomes at the day-ahead market ({\scshape s1,s2,s3,s4}). Furthermore, it can be observed that, while in the ConvMC model, the day-ahead market is cleared independently of the balancing market, StochMC determines the day-ahead dispatch accounting for the potential balancing cost of the system. 

Mathematically, the investment models depicted in Figure \ref{fig:FigInv2markets} are cast as hierarchical optimization problems. The upper-level problem determines the optimal expansion decisions of stochastic units to maximize the power producer expected profit. The lower-level problems represent the clearing of both the day-ahead and balancing markets for different conditions of the power system. Consequently, these models allow us to investigate two relevant aspects. Firstly, the impact of imbalance costs on generation expansion decisions of stochastic generating units can be assessed. This is accomplished by endogenously modeling, within the investment problem, the impact of capacity expansion decisions on the sequential clearing of the day-ahead and balancing markets. Secondly, these models can be used to evaluate the effect of the market-clearing mechanism on generating expansion decisions of power producers and to draw conclusions about which market design provides the appropriate incentives for investments in stochastic generating units. 

The mathematical formulation of the investment model considering a conventional market clearing is presented below. Relevant notation is provided in Table \ref{tab:Notation2}, where some of the symbols of formulation \eqref{Bilevel1Stage} are used here to refer to the day-ahead market. For example, while $\widehat{\bm{w}}_{s}$ in formulation \eqref{Bilevel1Stage} denotes the stochastic production in scenario $s$, here it represents the \emph{conditional expectation} of the stochastic production  used in ConvMC to clear the day-ahead market. We then denote by $\widetilde{\bm{w}}_{sr}$ the realized value of the stochastic production at the balancing stage.

\begin{table}[htbp]
	\centering \caption{Notation of optimization problem \eqref{B2S_Conv_upperOF}--\eqref{B2S_Conv_lowerlowerC3}}
		\begin{tabular}{c l}		
		\hline		
		\multicolumn{2}{l}{Numbers and sets}\\		
		\hspace{5mm} $s$   & Scenario index for day-ahead market conditions\\ 
		\hspace{5mm} $r$   & Scenario index for balancing market conditions\\	
		\multicolumn{2}{l}{Parameters}\\
		\hspace{5mm}$\bm{c}$ & $n_B\times1$ vector of bid and offer prices in the day-ahead market (\$/MWh) \\	
		\hspace{5mm}$\bm{c^{+}}$ & $n_B\times1$ vector of price offers for energy sale in the balancing market (\$/MWh)\\	
		\hspace{5mm}$\bm{c^{-}}$ & $n_B\times1$ vector of price offers for energy purchase in the balancing market (\$/MWh)\\	
		\hspace{5mm}$\widebar{\bm{f}_{sr}}$ & $n_L\times1$ vector of line capacities in balancing scenario $sr$ (MW) \\	
		\hspace{5mm}$\underline{\bm{p}_{sr}}$ & $n_B\times1$ vector of minimum dispatch limits in balancing scenario $sr$ (MW) \\
		\hspace{5mm}$\widebar{\bm{p}_{sr\hspace{1mm}}}$ & $n_B\times1$ vector of maximum dispatch limits in balancing scenario $sr$ (MW) \\
		\hspace{5mm}$\widehat{\bm{w}}_{s}$ & $n_B\times1$ vector of forecast stochastic power production in day-ahead scenario $s$ (p.u.) \\	
		\hspace{5mm}$\widetilde{\bm{w}}_{sr}$ & $n_B\times1$ vector of realized stochastic power production in balancing scenario $sr$ (p.u.) \\	
		\hspace{5mm}$\pi_s$ & Probability of occurrence of day-ahead scenario $s$ \\	
		\hspace{5mm}$\pi_{sr}$ & Probability of occurrence of balancing scenario $sr$ \\					
		\multicolumn{2}{l}{Variables}\\		
		\hspace{5mm}$\bm{p}_{s}$ & $n_B\times1$ vector of dispatched quantities in day-ahead scenario $s$ (MWh)\\
		\hspace{5mm}$\bm{p^{+}}_{sr}$ & $n_B\times1$ vector of increases in power injections in balancing scenario $sr$ (MWh)\\
		\hspace{5mm}$\bm{p^{-}}_{sr}$ & $n_B\times1$ vector of decreases in power injections in balancing scenario $sr$ (MWh)\\		
		\hspace{5mm}$\Delta\bm{p}^W_{sr}$ & $n_B\times1$ vector of dispatched stochastic power in balancing scenario $sr$ (MWh)\\
		\hspace{5mm}$\bm{\delta}_{s}$ & $n_B\times1$ vector of voltage angles in day-ahead scenario $s$ (rad)\\
		\hspace{5mm}$\Delta\bm{\delta}_{sr}$ & $n_B\times1$ vector of voltage angle differences in balancing scenario $sr$ (rad)\\	
		\hspace{5mm}$\bm{\lambda}_{s}$ & $n_B\times1$ vector of nodal prices in the day-ahead market scenario $s$ (\$/MWh)\\	
		\hspace{5mm}$\bm{\lambda}^R_{sr}$ & $n_B\times1$ vector of nodal prices in the balancing market in scenario $sr$ (\$/MWh)	\\
		\hline	
\end{tabular}\label{tab:Notation2}\end{table}

\begin{subequations}\label{Bilevel2StageConv}
\begin{align}
&\underset{\bm{U}}{{\rm Max}} \quad \sum_{s} \pi_{s} \left( \widehat{\bm{w}}^T_{s} diag(\bm{q}^T\bm{U})   \bm{\lambda}_{s} + \sum_r \left(\Delta\bm{p}^{W}_{sr}\right)^T \bm{\lambda}^R_{sr} \right) - \bm{i}^T\bm{U}\bm{1}  \label{B2S_Conv_upperOF} \hspace{55mm}\\
&{\rm s.t.} \quad \bm{U}\bm{1} \leq \bm{1} \label{B2S_Conv_upperC1} 
\end{align}
\vspace{-8mm}
\begin{empheq}[left={\hspace{7mm}\left(\Phi^D,\Phi^B \right) \in \text{arg}}\empheqlbrace,right=\empheqrbrace \forall s{,}]{align}
& \underset{\bm{p}_{s},\bm{\delta}_{s},\bm{p^{+}}_{sr},\bm{p^{-}}_{sr},\Delta\bm{p}^W_{sr},\Delta\bm{\delta}_{sr}}{{\rm Min}} \ \bm{c}^T\bm{p}_{s} + \sum_r \pi_{sr} \left( \bm{c}_+^{T}\bm{p^{+}}_{sr} + \bm{c}_-^{T}\bm{p^{-}}_{sr} \right) \label{B2S_Conv_lowerOF} \\
&  {\rm s.t.} \quad \underline{\bm{p}_{sr}} \leq \bm{p}_{s} + \bm{p^{+}}_{sr} - \bm{p^{-}}_{sr} \leq \widebar{\bm{p}_{sr}}: \underline{\bm{\xi}_{sr}}, \widebar{\bm{\xi}_{sr}}, \quad \forall r \label{B2S_Conv_lowerC1} \\
&  \hspace{8mm} -\widebar{\bm{f}_{sr}} \leq \bm{B}\bm{A} \left( \bm{\delta}_{s} + \Delta\bm{\delta}_{sr} \right) \leq \widebar{\bm{f}_{sr}}: \underline{\bm{\mu}_{sr}},\widebar{\bm{\mu}_{sr}} , \quad \forall r  \label{B2S_Conv_lowerC2} \\
&  \hspace{8mm} 0 \leq diag(\bm{q}^T\bm{U})\widehat{\bm{w}}_{s} + \Delta\bm{p}^W_{sr} \leq diag(\bm{q}^T\bm{U})\widetilde{\bm{w}}_{sr}: \underline{\bm{\psi}_{sr}}, \widebar{\bm{\psi}_{sr}}, \quad \forall r  \label{B2S_Conv_lowerC4} \\
&  \hspace{8mm} \bm{p^{+}}_{sr} - \bm{p^{-}}_{sr} +  \Delta\bm{p}^W_{sr} + \bm{A}^T\bm{B}\bm{A}\Delta\bm{\delta}_{sr} = \bm{0}: \bm{\lambda}^R_{sr}, \quad \forall r  \label{B2S_Conv_lowerC3} \\
& \phantom{A} \nonumber \\
& \phantom{A} \nonumber \\
& \phantom{A} \nonumber \\
& \phantom{A} \nonumber \\
& \phantom{A} \nonumber \\
& \phantom{A} \nonumber 
\end{empheq}
\vspace{-32mm}
\begin{empheq}[left={ \hspace{22mm}\Phi^D \in \text{arg}}\empheqlbrace,right=\empheqrbrace]{align}
\underset{\bm{p}_{s},\bm{\delta}_{s}}{{\rm  Min}} \quad & \bm{c}^T\bm{p}_{s}  \label{B2S_Conv_lowerlowerOF} \\
{\rm s.t.} \quad & \underline{\bm{p}_s} \leq \bm{p}_{s} \leq \widebar{\bm{p}_{s\hspace{1mm}}}: \underline{\bm{\alpha}_s}, \widebar{\bm{\alpha}_s} \label{B2S_Conv_lowerlowerC1} \\
& -\widebar{\bm{f}_s} \leq \bm{B}\bm{A}\bm{\delta}_{s} \leq \widebar{\bm{f}_s}: \underline{\bm{\beta}_s}, \widebar{\bm{\beta}_s} \label{B2S_Conv_lowerlowerC2} \\
& \bm{p}_s +   diag(\bm{q}^T\bm{U})\widehat{\bm{w}}_{s} + \bm{A}^T\bm{B}\bm{A}\bm{\delta}_{s} = \bm{0}: \bm{\lambda}_s \label{B2S_Conv_lowerlowerC3} 
\end{empheq}
\end{subequations}

\vspace{5mm}

\noindent where $\Phi^D = \{ \bm{\lambda}_{s},\bm{p}_s,\bm{\delta}_s \}$ and $\Phi^B = \{ \bm{\lambda}^R_{sr},\bm{p^{+}}_{sr},\bm{p^{-}}_{sr},\Delta\bm{p}^W_{sr},\Delta\bm{\delta}_{sr}\}$.

Model \eqref{Bilevel2StageConv} is a tri-level programming problem that determines the optimal investment decisions of a stochastic power producer assuming that the day-ahead and balancing markets are cleared independently.
The upper-level \eqref{B2S_Conv_upperOF}--\eqref{B2S_Conv_upperC1} optimizes the generation expansion plan $\bm{U}$ to maximize the producer's expected profit, which now includes the revenue from trading in both the day-ahead and balancing markets. Note that, as in \eqref{Bilevel1Stage}, constraint \eqref{B2S_Conv_upperC1} enforces that each wind farm can be located at only one node in the system. 

The second-level optimization problems \eqref{B2S_Conv_lowerOF}--\eqref{B2S_Conv_lowerlowerC3} represent the market-clearing formulation in which the day-ahead plus the expected balancing costs are minimized for each scenario $s$. Note that the outcome of the balancing market (re-dispatch of power producers and the balancing price) for a particular day-ahead scenario $s$ is influenced by the realization of a new set of random parameters (forecast errors of stochastic production and demand level) that unfold at the balancing stage. Constraint \eqref{B2S_Conv_lowerC1} defines
the feasible set of balancing adjustments for conventional units and demand. Likewise, transmission capacity limits are enforced by constraint \eqref{B2S_Conv_lowerC2}. Constraint \eqref{B2S_Conv_lowerC4} ensures that the sum of the stochastic power dispatched in the day-ahead market plus that cleared in the balancing market is positive and never higher than the actual power production for each scenario. Constraint \eqref{B2S_Conv_lowerC3} represents the nodal power balance equations at the balancing stage. The price at which imbalances are settled is computed as the probability-removed dual variable of this constraint \citep{wong2007pricing}. 

Lastly, the third level \eqref{B2S_Conv_lowerlowerOF}--\eqref{B2S_Conv_lowerlowerC3} represents the clearing of the day-ahead market exactly as in \eqref{Bilevel1Stage}. Observe that by including this optimization problem as a set of constraints in the overall market-clearing problem \eqref{B2S_Conv_lowerOF}--\eqref{B2S_Conv_lowerlowerC3}, we ensure that the resulting day-ahead quantities and prices ($\bm{p}_s$ and $\bm{\lambda}_s$) are those that minimize the day-ahead scheduling costs for each scenario $s$ (which, in general, do not coincide with those that minimize the expected total system cost). This formulation captures then the sequential and independent clearing of the day-ahead and balancing markets and provides a dispatch equivalent to that of a market design in which reserve capacity is auctioned after the closure of the day-ahead energy market. Note that, following current practice, the dispatched stochastic power in the day-ahead market is assumed to be equal to the forecast production $\widehat{\bm{w}}_{s}$ \citep{bouffard2008stochastic}. It is also worth mentioning that, although $\bm{p}_{s}$ is included in the decision variable set of objective function \eqref{B2S_Conv_lowerOF} for clarity, optimization problem \eqref{B2S_Conv_lowerlowerOF}--\eqref{B2S_Conv_lowerlowerC3} normally provides a unique solution for day-ahead dispatch quantities.

Before  presenting the procedure to solve \eqref{Bilevel2StageConv}, some comments are in order. First, it should be stressed that we assume an energy-only market and therefore, capacity payments are not considered here. Second, the generating expansion models we work with disregard non-convexities of generating units such as capacity limitations and large start-up costs. 

In order to solve the tri-level optimization problem \eqref{Bilevel2StageConv}, we need first to acknowledge that for a fix value of the upper-level decision variable $\bm{U}$, optimization problem \eqref{B2S_Conv_lowerlowerOF}--\eqref{B2S_Conv_lowerlowerC3} is linear and therefore, it can be replaced with its primal constraints, the constraints of its corresponding dual optimization problem, and the primal-dual strong duality theorem \citep{motto2005mixed} as follows:
\begin{subequations}\label{Bilevel2StageConvPD1}
\begin{align}
&\underset{\bm{U}}{{\rm Max}} \quad \sum_{s} \pi_{s} \left( \widehat{\bm{w}}^T_{s}diag(\bm{q}^T\bm{U}) \bm{\lambda}_{s} + \sum_r \left(\Delta\bm{p}^W_{sr}\right)^T\bm{\lambda}^R_{sr} \right) - \bm{i}^T\bm{U}\bm{1}  \label{B2S_Conv_PD1_upperOF} \hspace{50mm}\\
&{\rm s.t.} \quad \bm{U}\bm{1} \leq \bm{1} \label{B2S_Conv_PD1_upperC1} 
\end{align}
\vspace{-8mm}
\begin{empheq}[left={\hspace{0mm}\left(\Phi^D,\Phi^B \right) \in \text{arg}}\empheqlbrace,right=\empheqrbrace \forall s{,}]{align}
& \underset{\bm{p}_{s},\bm{\delta}_{s},\bm{p^{+}}_{sr},\bm{p^{-}}_{sr},\Delta\bm{p}^W_{sr},\Delta\bm{\delta}_{sr}}{{\rm Min}} \ \bm{c}^T\bm{p}_{s} + \sum_r \pi_{sr} \left( \bm{c}_+^{T}\bm{p^{+}}_{sr} + \bm{c}_-^{T}\bm{p^{-}}_{sr} \right) \label{B2S_Conv_PD1_lowerOF} \\
&  {\rm s.t.} \hspace{2mm} {\rm (\ref{B2S_Conv_lowerC1})-(\ref{B2S_Conv_lowerC4})} \nonumber \\
&  \hspace{8mm} \underline{\bm{p}_s} \leq \bm{p}_{s} \leq \widebar{\bm{p}_{s\hspace{1mm}}}: \widecheck{\bm{\alpha}}_s, \widehat{\bm{\alpha}}_s \label{B2S_Conv_PD1_lowerlowerC1} \\
&  \hspace{8mm} -\widebar{\bm{f}_s} \leq \bm{B}\bm{A}\bm{\delta}_{s} \leq \widebar{\bm{f}_s}: \widecheck{\bm{\beta}}_s, \widehat{\bm{\beta}}_s \label{B2S_Conv_PD1_lowerlowerC2} \\
& \hspace{8mm} \bm{p}_s +   diag(\bm{q}^T\bm{U})\widehat{\bm{w}}_{s} + \bm{A}^T\bm{B}\bm{A}\bm{\delta}_{s} = \bm{0}: \widehat{\bm{\lambda}}_s  \label{B2S_Conv_PD1_lowerlowerC3} \\
& \hspace{8mm} \underline{\bm{\alpha}_s} + \widebar{\bm{\alpha}_s} + \bm{\lambda}_s = \bm{c}: \bm{\eta}_s  \label{B2S_Conv_PD1_lowerlowerD1} \\
& \hspace{8mm} \bm{A}^T\bm{B} \left( \underline{\bm{\beta}_s} + \widebar{\bm{\beta}_s} \right) + \bm{A}^T\bm{B}\bm{A}\bm{\lambda}_s = 0 : \bm{\tau}_s \label{B2S_Conv_PD1_lowerlowerD2} \\
& \hspace{8mm} \bm{c}^T\bm{p}_{s} =  \underline{\bm{p}_s}^T\underline{\bm{\alpha}_s} + \widebar{\bm{p}_{s\hspace{1mm}}}^T\widebar{\bm{\alpha}_s} -
\widebar{\bm{f}_s}^T\underline{\bm{\beta}_s} +
\widebar{\bm{f}_s}^T\widebar{\bm{\beta}_s}  - \widehat{\bm{w}}^T_{s} diag(\bm{q}^T\bm{U})\bm{\lambda}_s: \phi_s \label{B2S_Conv_PD1_lowerlowerD3} \\
& \hspace{8mm} \underline{\bm{\alpha}_s}, \underline{\bm{\beta}_s} \geq \bm{0} \label{B2S_Conv_PD1_lowerlowerD4} \\
& \hspace{8mm} \widebar{\bm{\alpha}_s}, \widebar{\bm{\beta}_s} \leq \bm{0}  \label{B2S_Conv_PD1_lowerlowerD5}
\end{empheq}
\end{subequations}

\vspace{5mm}

\noindent where constraints  \eqref{B2S_Conv_PD1_lowerlowerC1}--\eqref{B2S_Conv_PD1_lowerlowerC3} are the primal constraints of the third-level problem; equations \eqref{B2S_Conv_PD1_lowerlowerD1}--\eqref{B2S_Conv_PD1_lowerlowerD2} represent the dual constraints corresponding to variables $\bm{p}_{s}$ and $\bm{\delta}_s$, respectively; equation \eqref{B2S_Conv_PD1_lowerlowerD3} ensures that the primal and dual formulation of the lower-level problems reach the same objective function at the optimal solution; and \eqref{B2S_Conv_PD1_lowerlowerD4} and \eqref{B2S_Conv_PD1_lowerlowerD5} are positive and negative variable declarations, respectively. Note also that new dual variables corresponding to constraints \eqref{B2S_Conv_PD1_lowerlowerC1}--\eqref{B2S_Conv_PD1_lowerlowerD3} have been included in \eqref{Bilevel2StageConvPD1}. For example, while $\widebar{\bm{\alpha}_s}$ in \eqref{Bilevel2StageConv} represents the sensitivity of the objective function \eqref{B2S_Conv_lowerlowerOF} to variations of $\widebar{\bm{p}_{s\hspace{1mm}}}$, $\widehat{\bm{\alpha}}_s$ in \eqref{Bilevel2StageConvPD1} quantifies the impact of marginal changes of $\widebar{\bm{p}_{s\hspace{1mm}}}$ on objective function \eqref{B2S_Conv_PD1_lowerOF}. 

Observe that, after replacing the third-level problem with \eqref{B2S_Conv_PD1_lowerlowerC1}--\eqref{B2S_Conv_PD1_lowerlowerD5}, optimization problem \eqref{B2S_Conv_PD1_lowerOF}--\eqref{B2S_Conv_PD1_lowerlowerD5} remains linear and therefore, it can be also replaced with its
primal constraints, the constraints of the corresponding dual problem, as well as the strong duality condition. It is worth clarifying that the dual variables of the original day-ahead market-clearing problem $\widebar{\bm{\alpha}_s}$, $\underline{\bm{\alpha}_s}$, $\widebar{\bm{\beta}_s}$, $\underline{\bm{\beta}_s}$, and $\bm{\lambda}_s$ turn into primal variables within  formulation \eqref{B2S_Conv_PD1_lowerOF}--\eqref{B2S_Conv_PD1_lowerlowerD5}. After replacing the second-level optimization problem with its primal and dual constraints plus the strong duality condition, the original three-level investment model \eqref{Bilevel2StageConv} is recast as an equivalent one-level optimization problem as presented below: 
\begin{subequations}\label{Bilevel2StageConvPD2}
\begin{align}
&\underset{\scriptsize \shortstack{$\bm{U},\Phi^D,\Phi^B$}}{{\rm Max}} \quad \sum_{s} \pi_{s} \left( \widehat{\bm{w}}^T_{s} diag(\bm{q}^T\bm{U})  \bm{\lambda}_{s} + \sum_r \left(\Delta\bm{p}^W_{sr}\right)^T \bm{\lambda}^R_{sr} \right) - \bm{i}^T\bm{U}\bm{1}  \label{B2S_Conv_PD2_upperOF}\\
&{\rm s.t.} \quad \bm{U}\bm{1} \leq \bm{1} \label{B2S_Conv_PD2_upperC1} \\
& \qquad \eqref{B2S_Conv_lowerC1}-\eqref{B2S_Conv_lowerC4}, \eqref{B2S_Conv_PD1_lowerlowerC1}-\eqref{B2S_Conv_PD1_lowerlowerD5} \nonumber \\
& \qquad  \widecheck{\bm{\alpha}}_s + \widehat{\bm{\alpha}}_s + \widehat{\bm{\lambda}}_s + \sum_r \left( \underline{\bm{\xi}_{sr}} + \widebar{\bm{\xi}_{sr}} \right) + \phi_s\bm{c} = \bm{c}, \quad \forall s \label{B2S_Conv_PD2_lowerD1} \\
& \qquad  \bm{A}^T\bm{B} \left( \widecheck{\bm{\beta}}_s + \widehat{\bm{\beta}}_s \right) + \bm{A}^T\bm{B}\bm{A}\widehat{\bm{\lambda}}_s + \sum_r\bm{A}^T\bm{B}\left( \underline{\bm{\mu}_{sr}} +\widebar{\bm{\mu}_{sr}} \right) = \bm{0}, \quad \forall s \label{B2S_Conv_PD2_lowerD2} \\
& \qquad  \underline{\bm{\xi}_{sr}} + \widebar{\bm{\xi}_{sr}}  + \bm{\lambda}^R_{sr} = \pi_{sr}\bm{c}_{+}, \quad \forall s, \forall r\label{B2S_Conv_PD2_lowerD3} \\
& \qquad  -\underline{\bm{\xi}_{sr}} - \widebar{\bm{\xi}_{sr}} - \bm{\lambda}^R_{sr} = \pi_{sr}\bm{c}_{-}, \quad \forall s, \forall r\label{B2S_Conv_PD2_lowerD4} \\
& \qquad  \bm{A}^T\bm{B}\left( \underline{\bm{\mu}_{sr}} +\widebar{\bm{\mu}_{sr}} \right) + \bm{A}^T\bm{B}\bm{A}\bm{\lambda}^R_{sr} = \bm{0}, \quad \forall s, \forall r\label{B2S_Conv_PD2_lowerD5} \\
& \qquad  \bm{\lambda}^R_{sr} + \underline{\bm{\psi}_{sr}} + \widebar{\bm{\psi}_{sr}} = \bm{0}, \quad \forall s, \forall r \label{B2S_Conv_PD2_lowerD6} \\
& \qquad  \bm{\eta}_s - \underline{\bm{p}_s}\phi_s = \bm{0}, \quad \forall s \label{B2S_Conv_PD2_lowerD7} \\
& \qquad  \bm{\eta}_s - \widebar{\bm{p}_{s\hspace{1mm}}}\phi_s = \bm{0}, \quad \forall s \label{B2S_Conv_PD2_lowerD8} \\
& \qquad  \bm{B}\bm{A}\bm{\tau}_s + \widebar{\bm{f}_s}\phi_s = \bm{0}, \quad \forall s \label{B2S_Conv_PD2_lowerD9} \\
& \qquad  \bm{B}\bm{A}\bm{\tau}_s - \widebar{\bm{f}_s}\phi_s = \bm{0}, \quad \forall s \label{B2S_Conv_PD2_lowerD10} \\
& \qquad \bm{\eta}_s + \bm{A}^T\bm{B}\bm{A}\bm{\tau}_s + \widehat{\bm{w}}^T_{s} diag(\bm{q}^T\bm{U})\phi_s = \bm{0}, \quad \forall s \label{B2S_Conv_PD2_lowerD11} \\
& \qquad  \bm{c}^T\bm{p}_{s} + \sum_r \pi_{sr} \left( \bm{c}_+^{T}\bm{p^{+}}_{sr} - \bm{c}_-^{T}\bm{p^{-}}_{sr} \right) = 
\widebar{\bm{p}_{s\hspace{1mm}}}^T\widehat{\bm{\alpha}}_s + 
\underline{\bm{p}_s}^T\widecheck{\bm{\alpha}}_s +
\widebar{\bm{f}_s}^T\widehat{\bm{\beta}}_s - 
\underline{\bm{f}_s}^T\widecheck{\bm{\beta}}_s -
\widehat{\bm{w}}^T_{s}diag(\bm{q}^T\bm{U})\widehat{\bm{\lambda}}_s +
\bm{c}^T\bm{\eta}_s + \nonumber\\
& \qquad + \sum_r \left( \widebar{\bm{p}_{sr}}^T\widebar{\bm{\xi}_{sr}} + 
\underline{\bm{p}_{sr}}^T\underline{\bm{\xi}_{sr}} + 
\widebar{\bm{f}_{sr}}^T\widebar{\bm{\mu}_{sr}} -
\underline{\bm{f}_{sr}}^T\underline{\bm{\mu}_{sr}} + diag(\bm{q}^T\bm{U})\left( \left(\widetilde{\bm{w}}_{sr}^T - \widehat{\bm{w}}^T_s \right)\widebar{\bm{\psi}_{sr}} - \widehat{\bm{w}}^T_s\underline{\bm{\psi}_{sr}} \right) \right), \quad \forall s \label{B2S_Conv_PD2_lowerD12} \\
&\qquad\widecheck{\bm{\alpha}}_s, \widecheck{\bm{\beta}}_s, \underline{\bm{\xi}_{sr}}, \underline{\bm{\mu}_{sr}}, \underline{\bm{\psi}_{sr}}  \geq \bm{0}, \quad \forall s \label{B2S_Conv_PD2_lowerD13} \\
&\qquad\widehat{\bm{\alpha}}_s, \widehat{\bm{\beta}}_s, \widebar{\bm{\xi}_{sr}}, \widebar{\bm{\mu}_{sr}}, \widebar{\bm{\psi}_{sr}} \leq \bm{0}, \quad \forall s.  \label{B2S_Conv_PD2_lowerD14}
\end{align}
\end{subequations}

Objective function \eqref{B2S_Conv_PD2_upperOF} maximizes the expected profit of newly installed stochastic generating units. Set of constraints \eqref{B2S_Conv_lowerC1}-\eqref{B2S_Conv_lowerC4} are primal constraints corresponding to the balancing market-clearing problem, while equations \eqref{B2S_Conv_PD1_lowerlowerC1}-\eqref{B2S_Conv_PD1_lowerlowerD5} represent the primal and dual constraints, together with the strong duality condition of the day-ahead market-clearing problem. Likewise,  \eqref{B2S_Conv_PD2_lowerD1}--\eqref{B2S_Conv_PD2_lowerD11} are the dual constraints associated with variables $\bm{p}_s$, $\bm{\delta}_s$, $\bm{p^{+}}_{sr}$, $\bm{p^{-}}_{sr}$, $\Delta\bm{\delta}_{sr}$, $\Delta\bm{p}^W_{sr}$, $\widebar{\bm{\alpha}_s}$, $\underline{\bm{\alpha}_s}$, $\widebar{\bm{\beta}_s}$, $\underline{\bm{\beta}_s}$, and $\bm{\lambda}_s$, in that order. Equation \eqref{B2S_Conv_PD2_lowerD12} imposes the strong duality condition related to optimization problem \eqref{B2S_Conv_PD1_lowerOF}--\eqref{B2S_Conv_PD1_lowerlowerD5}. Finally, \eqref{B2S_Conv_PD2_lowerD13} and \eqref{B2S_Conv_PD2_lowerD14} are positive and negative variable declarations.

Observe that, besides the second term of the objective function, the rest of constraints in \eqref{Bilevel2StageConvPD2} contain either linear terms or products of binary and continuous variables, which can be linearized as described in \citet{floudas1995nonlinear}. In order to linearize the term $\left(\Delta\bm{p}^W_{sr}\right)^T\bm{\lambda}^R_{sr}$ in objective function \eqref{B2S_Conv_PD2_upperOF}, we first multiply equation \eqref{B2S_Conv_PD2_lowerD6} by $\left(\Delta\bm{p}^W_{sr}\right)^T$, thus obtaining
\begin{align}
& \left(\Delta\bm{p}^W_{sr}\right)^T \bm{\lambda}^R_{sr} = - \left(\Delta\bm{p}^W_{sr}\right)^T\left(\underline{\bm{\psi}_{sr}} + \widebar{\bm{\psi}_{sr}} \right).  \label{Lin2_1}
\end{align}
Furthermore, the complementarity conditions of \eqref{B2S_Conv_lowerC4} imply that 
\begin{align}
& - \left(\Delta\bm{p}^W_{sr}\right)^T\widebar{\bm{\psi}_{sr}} =  \left( \widehat{\bm{w}}_{s}- \widetilde{\bm{w}}_{sr} \right)^Tdiag(\bm{q}^T\bm{U})\widebar{\bm{\psi}_{sr}}, \label{Lin2_2} \\
& - \left(\Delta\bm{p}^W_{sr}\right)^T\underline{\bm{\psi}_{sr}} =  \widehat{\bm{w}}^T_{s}diag(\bm{q}^T\bm{U})\underline{\bm{\psi}_{sr}}. \label{Lin2_3}
\end{align} 
Replacing \eqref{Lin2_2} and \eqref{Lin2_3} in \eqref{Lin2_1}, we obtain the following equivalent linear expression of the nonlinear term included in  objective function \eqref{B2S_Conv_PD2_upperOF}
\begin{align}
&  \left(\Delta\bm{p}^W_{sr}\right)^T \bm{\lambda}^R_{sr} =  \widehat{\bm{w}}^T_{s}diag(\bm{q}^T\bm{U})\big( \widebar{\bm{\psi}_{sr}} +  \underline{\bm{\psi}_{sr}} \big) - \widetilde{\bm{w}}^T_{sr}diag(\bm{q}^T\bm{U})\widebar{\bm{\psi}_{sr}}.
\label{Lin2_4}
\end{align}
By using \eqref{B2S_Conv_PD2_lowerD6} again, we can rewrite objective function \eqref{B2S_Conv_PD2_upperOF} as
\begin{align}
&  \sum_{s} \pi_{s}\widehat{\bm{w}}^T_{s}diag(\bm{q}^T\bm{U})   \Big(\bm{\lambda}_{s}-\sum_r \bm{\lambda}^R_{sr} \Big) 
- \sum_{sr} \pi_{s}\widetilde{\bm{w}}_{sr}^Tdiag(\bm{q}^T\bm{U})\widebar{\bm{\psi}_{sr}} - \bm{i}^T\bm{U}\bm{1}, \label{Lin2_5}
\end{align}
\noindent which only contains products of binary and continuous variables. Therefore, optimization problem \eqref{Bilevel2StageConvPD2} is now formulated as a mixed-integer linear programming problem that can be solved using commercially available optimization software. 

Alternatively, we formulate next the investment optimization problem of a stochastic power producer participating in a market in which day-ahead dispatch decisions are determined according to a stochastic market-clearing procedure \citep{pritchard2010single, Morales2012}: 
\begin{subequations}\label{Bilevel2StageSto}
\begin{align}
&\underset{\bm{U}}{{\rm Max}} \quad \sum_{s} \pi_{s} \left( \left(\bm{p}^W_{s}\right)^T  \bm{\lambda}_{s} + \sum_r \left(\Delta\bm{p}^W_{sr}\right)^T \bm{\lambda}^R_{sr} \right) - \bm{i}^T\bm{U}\bm{1}  \label{B2S_Sto_upperOF} \hspace{55mm}\\
&{\rm s.t.} \quad \bm{U}\bm{1} \leq \bm{1} \label{B2S_Sto_upperC1} 
\end{align}
\vspace{-8mm}
\begin{empheq}[left={\hspace{10mm}\Phi^{DB}  \in \text{arg}}\empheqlbrace,right=\empheqrbrace \forall s{,}]{align}
& \underset{\bm{p}_{s},\bm{p}^W_{s},\bm{\delta}_{s},\bm{p^{+}}_{sr},\bm{p^{-}}_{sr},\Delta\bm{p}^W_{sr},\Delta\bm{\delta}_{sr}}{{\rm Min}} \ \bm{c}^T\bm{p}_{s} + \sum_r \pi_{sr} \left( \bm{c}_+^{T}\bm{p^{+}}_{sr} + \bm{c}_-^{T}\bm{p^{-}}_{sr} \right) \label{B2S_Sto_lowerOF} \\
& \hspace{10mm} {\rm s.t.} \quad \underline{\bm{p}_s} \leq \bm{p}_{s} \leq \widebar{\bm{p}_{s\hspace{1mm}}}: \underline{\bm{\alpha}_s}, \widebar{\bm{\alpha}_s} \label{B2S_Sto_lowerC1} \\
& \hspace{19mm} -\widebar{\bm{f}_s} \leq \bm{B}\bm{A}\bm{\delta}_{s} \leq \widebar{\bm{f}_s}: \underline{\bm{\beta}_s}, \widebar{\bm{\beta}_s} \label{B2S_Sto_lowerC2} \\
& \hspace{19mm} \bm{p}_s +   \bm{p}^W_{s} + \bm{A}^T\bm{B}\bm{A}\bm{\delta}_{s} = \bm{0}: \bm{\lambda}_s \label{B2S_Sto_lowerC3} \\
& \hspace{19mm} \underline{\bm{p}_{sr}} \leq \bm{p}_{s} + \bm{p^{+}}_{sr} - \bm{p^{-}}_{sr} \leq \widebar{\bm{p}_{sr}}: \underline{\bm{\xi}_{sr}}, \widebar{\bm{\xi}_{sr}}, \quad \forall r \label{B2S_Sto_lowerC4} \\
& \hspace{19mm} -\widebar{\bm{f}_{sr}} \leq \bm{B}\bm{A} \left( \bm{\delta}_{s} + \Delta\bm{\delta}_{sr} \right) \leq \widebar{\bm{f}_{sr}}: \underline{\bm{\mu}_{sr}},\widebar{\bm{\mu}_{sr}} , \quad \forall r \label{B2S_Sto_lowerC5} \\
& \hspace{19mm} \bm{0} \leq \bm{p}^W_{s} + \Delta\bm{p}^W_{sr} \leq diag(\bm{q}^T\bm{U})\bm{\widetilde{w}_{sr}}: \underline{\bm{\psi}_{sr}}, \widebar{\bm{\psi}_{sr}} , \quad \forall r  \label{B2S_Sto_lowerC6} \\
& \hspace{19mm} \bm{p^{+}}_{sr} - \bm{p^{-}}_{sr} +  \Delta\bm{p}^W_{sr} + \bm{A}^T\bm{B}\bm{A}\Delta\bm{\delta}_{sr} = \bm{0}: \bm{\lambda}^R_{sr}, \quad \forall r  \label{B2S_Sto_lowerC7}
\end{empheq}
\end{subequations}

\vspace{5mm}

\noindent \noindent where $\bm{p}^W_{s}$ corresponds to the day-ahead dispatch of the stochastic producers in scenario $s$ and $\Phi^{DB} = \{ \bm{\lambda}_{s},\bm{p}_s,\bm{p}^W_{s}, \bm{\delta}_s,\bm{\lambda}^R_{sr},\bm{p^{+}}_{sr},\bm{p^{-}}_{sr},\Delta\bm{p}^W_{sr},\Delta\bm{\delta}_{sr}\}$. Observe that the two levels \eqref{B2S_Conv_lowerOF}--\eqref{B2S_Conv_lowerC3} and \eqref{B2S_Conv_lowerlowerOF}--\eqref{B2S_Conv_lowerlowerC3} of the conventional market-clearing problem considered in investment model \eqref{Bilevel2StageConv} have been merged into the single-level problem \eqref{B2S_Sto_lowerOF}--\eqref{B2S_Sto_lowerC7}, which jointly minimizes the day-ahead dispatch cost and the expected cost of balancing the system. Problem \eqref{B2S_Sto_lowerOF}--\eqref{B2S_Sto_lowerC7} is actually a two-stage stochastic programming problem, where the day-ahead dispatch decisions for each scenario $s$ ($\bm{p}_s,\bm{p}^W_{s},\bm{\delta}_{s}$) are computed so that the subsequent re-dispatch actions for each balancing scenario $r$ ($\bm{p^{+}}_{sr},\bm{p^{-}}_{sr},\Delta\bm{ w}^T_{sr},\Delta\bm{\delta}_{sr}$) minimize the expected total system operation cost as expressed in objective function \eqref{B2S_Sto_lowerOF}. Consequently, the clearing mechanism \eqref{B2S_Sto_lowerOF}--\eqref{B2S_Sto_lowerC7} provides a generation dispatch equivalent to that of a market design where reserve capacity is determined and allocated in concurrence with the day-ahead energy dispatch using a model that weighs the predicted balancing cost for each realization of the stochastic power production \citep{bouffard2008stochastic,Morales09b}. Note also that reserve capacity bids are disregarded here since they do not represent an intrinsic cost for the generators \citep{papavasiliou2011reserve}. 

Similarly to the previous model, the upper-level objective function \eqref{B2S_Sto_upperOF} maximizes the expected revenue of a power producer investing in new stochastic generating units, the first and second term corresponding to the revenue obtained in the day-ahead and balancing market, respectively. Constraint \eqref{B2S_Sto_upperC1} states that each generation expansion project can be only placed at one node of the network.

The lower-level optimization problem \eqref{B2S_Sto_lowerOF}--\eqref{B2S_Sto_lowerC7} represents a day-ahead market-clearing model that accounts for the impact of day-ahead dispatch decisions on the future balancing costs of the system. The set of constraints \eqref{B2S_Sto_lowerC1}--\eqref{B2S_Sto_lowerC2} and \eqref{B2S_Sto_lowerC4}--\eqref{B2S_Sto_lowerC6} impose limits on day-ahead and balancing dispatch decisions, respectively. Similarly, the power balance at each scenario is enforced at the day-ahead and balancing stages through equations \eqref{B2S_Sto_lowerC3} and \eqref{B2S_Sto_lowerC7}, respectively, being $\bm{\lambda}_s$ and $\bm{\lambda}^R_{sr}$ the resulting clearing electricity prices at each market floor. 

It should be noticed that, unlike the market-clearing model \eqref{B2S_Conv_lowerOF}--\eqref{B2S_Conv_lowerlowerC3}, in which the day-ahead dispatch of the stochastic generating units is assumed to be equal to their expected power production $\widehat{\bm{w}}_{s}$, market clearing \eqref{B2S_Sto_lowerOF}--\eqref{B2S_Sto_lowerC7} can dispatch stochastic units to values different from their expected production through the continuous variable $\bm{p}^W_{s}$. In doing so, the day-ahead dispatch of stochastic units can be adapted to the availability and cost of balancing resources provided by flexible generating units. Besides, observe that, while model \eqref{Bilevel2StageConv} provides the same day-ahead dispatch decisions regardless of how significant the forecast errors are or how much flexibility can be provided by the conventional generating units at the balancing stage, the market-clearing model \eqref{B2S_Sto_lowerOF}--\eqref{B2S_Sto_lowerC7} can dispatch flexible conventional units out of merit in the day-ahead market to efficiently cope with energy imbalances during the real-time operation of the system. StochMC results consequently in less volatile prices at the balancing stage than ConvMC, thus reducing the balancing costs incurred by stochastic power producers due to their production forecast errors. This has, in turn, an impact on the investment decisions made by stochastic power producers. 

In order to solve the bi-level optimization problem \eqref{Bilevel2StageSto}, we construct an equivalent one-level optimization problem. To this end, and provided that model \eqref{B2S_Sto_lowerOF}--\eqref{B2S_Sto_lowerC7} is linear for given investment decisions $\bm{U}$, we can replace the set of lower-level problems with their primal and dual constraints plus the strong duality conditions as follows:

\begin{subequations}\label{Bilevel2StageStoPD}
\begin{align}
&\underset{\scriptsize \shortstack{$\bm{U},\Phi^{DB}$}}{{\rm Maximize}} \quad \sum_{s} \pi_{s} \left( \left(\bm{p}^W_{s}\right)^T  \bm{\lambda}_{s} + \sum_r \left(\Delta\bm{p}^W_{sr}\right)^T \bm{\lambda}^R_{sr} \right) - \bm{i}^T\bm{U}\bm{1}  \label{B2S_StoPD_upperOF}\\
&{\rm s.t.} \quad \bm{U}\bm{1} \leq \bm{1} \label{B2S_StoPD_upperC1} \\
& \qquad \eqref{B2S_Sto_lowerC1}-\eqref{B2S_Sto_lowerC7} \nonumber \\
& \qquad \underline{\bm{\alpha}_s} + \widebar{\bm{\alpha}_s} + \bm{\lambda}_s + \sum_r \left( \underline{\bm{\xi}_{sr}} + \widebar{\bm{\xi}_{sr}} \right) = \bm{c}, \quad \forall s \label{B2S_StoPD_lowerD1} \\
& \qquad \bm{A}^T\bm{B}\left( \underline{\bm{\beta}_s} + \widebar{\bm{\beta}_s} \right) + \bm{A}^T\bm{B}\bm{A}\bm{\lambda}_s + \sum_r\bm{A}^T\bm{B}\left( \underline{\bm{\mu}_{sr}} + \widebar{\bm{\mu}_{sr}} \right) = 0, \quad \forall s \label{B2S_StoPD_lowerD2} \\
& \qquad \bm{\lambda}_s + \sum_r \left( \underline{\bm{\psi}_{sr}} + \widebar{\bm{\psi}_{sr}} \right) = 0, \quad \forall s \label{B2S_StoPD_lowerD3} \\
& \qquad \underline{\bm{\xi}_{sr}} + \widebar{\bm{\xi}_{sr}} + \bm{\lambda}^R_{sr} = \pi_{sr}\bm{c}_{+}, \quad \forall s, \forall r\label{B2S_StoPD_lowerD4} \\
& \qquad - \underline{\bm{\xi}_{sr}} - \widebar{\bm{\xi}_{sr}} - \bm{\lambda}^R_{sr} = \pi_{sr}\bm{c}_{-}, \quad \forall s, \forall r\label{B2S_StoPD_lowerD5} \\
& \qquad \bm{A}^T\bm{B} \left( \underline{\bm{\mu}_{sr}} + \widebar{\bm{\mu}_{sr}} \right) + \bm{A}^T\bm{B}\bm{A}\bm{\lambda}^R_{sr} = 0, \quad \forall s, \forall r\label{B2S_StoPD_lowerD6} \\
& \qquad \underline{\bm{\psi}_{sr}} + \widebar{\bm{\psi}_{sr}} + \bm{\lambda}^R_{sr} = 0, \quad \forall s, \forall r \label{B2S_StoPD_lowerD7} \\
& \qquad \bm{c}^T\bm{p}_{s} + \sum_r \pi_{sr} \left( \bm{c}_+^{T}\bm{p^{+}}_{sr} + \bm{c}_-^{T}\bm{p^{-}}_{sr} \right) = 
\widebar{\bm{p}_{s\hspace{1mm}}}^T\widebar{\bm{\alpha}_s} + 
\underline{\bm{p}_s}^T\underline{\bm{\alpha}_s} +
\widebar{\bm{f}_s}^T\widebar{\bm{\beta}_s} - 
\underline{\bm{f}_s}^T\underline{\bm{\beta}_s} + \nonumber\\
& \qquad + \sum_r \left( \widebar{\bm{p}_{sr}}^T\widebar{\bm{\xi}_{sr}} + 
\underline{\bm{p}_{sr}}^T\underline{\bm{\xi}_{sr}} + 
\widebar{\bm{f}_{sr}}^T\widebar{\bm{\mu}_{sr}} -
\underline{\bm{f}_{sr}}^T\underline{\bm{\mu}_{sr}} + diag(\bm{q}^T\bm{U}) \widetilde{\bm{w}}_{sr}^T\widebar{\bm{\psi}_{sr}} \right), \quad \forall s \label{B2S_StoPD_lowerD8} \\
&\qquad\underline{\bm{\alpha}_s}, \underline{\bm{\beta}_s}, \underline{\bm{\xi}_{sr}}, \underline{\bm{\mu}_{sr}}, \underline{\bm{\psi}_{sr}}  \geq \bm{0}, \quad \forall s, \forall r \label{B2S_StoPD_lowerD9} \\
&\qquad\widebar{\bm{\alpha}_s}, \widebar{\bm{\beta}_s}, \widebar{\bm{\xi}_{sr}}, \widebar{\bm{\mu}_{sr}}, \widebar{\bm{\psi}_{sr}} \leq \bm{0}, \quad \forall s, \forall r. \label{B2S_StoPD_lowerD10}
\end{align}
\end{subequations}
where \eqref{B2S_StoPD_lowerD1}--\eqref{B2S_StoPD_lowerD7} are the dual constraints associated with primal variables $\bm{p}_s$, $\bm{\delta}_s$, $\bm{p}^W_{s}$, $\bm{p^{+}}_{sr}$, $\bm{p^{-}}_{sr}$, $\Delta\bm{\delta}_{sr}$, $\Delta\bm{p}^W_{sr}$, in that order; \eqref{B2S_StoPD_lowerD8} corresponds to the strong duality condition; and \eqref{B2S_StoPD_lowerD9} and \eqref{B2S_StoPD_lowerD10} are positive and negative variable declarations, respectively. Observe that all constraints in \eqref{Bilevel2StageStoPD} contain linear terms or products of binary and continuous variables that can be easily linearized. However, objective function \eqref{B2S_StoPD_upperOF} includes non-linear products that are linearized as follows. First, the complementarity conditions associated with constraints \eqref{B2S_Sto_lowerC6} imply that
\begin{align}
& \left(\bm{p}^W_{s}\right)^T\widebar{\bm{\psi}_{sr}} = -\left( \Delta\bm{p}^W_{sr} - diag(\bm{q}^T\bm{U})\widetilde{\bm{w}}_{sr} \right)^T\widebar{\bm{\psi}_{sr}}, \label{Lin1_1} \\
&  \left(\bm{p}^W_{s}\right)^T\underline{\bm{\psi}_{sr}} = - \left(\Delta\bm{p}^W_{sr}\right)^T\underline{\bm{\psi}_{sr}}. \label{Lin1_2}
\end{align} 
Multiplying the dual constraint \eqref{B2S_StoPD_lowerD3} by $\left(\bm{p}^W_{s}\right)^T$ and using \eqref{Lin1_1} and \eqref{Lin1_2} we obtain
\begin{align}
&  \left(\bm{p}^W_{s}\right)^T  \bm{\lambda}_{s} = \sum_r \left( \Delta\bm{p}^W_{sr} - diag(\bm{q}^T\bm{U})\widetilde{\bm{w}}_{sr} \right)^T\widebar{\bm{\psi}_{sr}} + \sum_r \left(\Delta\bm{p}^W_{sr}\right)^T\underline{\bm{\psi}_{sr}}.  \label{Lin1_3}
\end{align}
Similarly, we can multiply \eqref{B2S_StoPD_lowerD7} by $\left(\Delta\bm{p}^W_{sr}\right)^T$, sum over $r$ and obtain
\begin{align}
& \sum_r \left(\Delta\bm{p}^W_{sr}\right)^T \bm{\lambda}^R_{sr} = -\sum_r \left(\Delta\bm{p}^W_{sr}\right)^T \left(  \underline{\bm{\psi}_{sr}} + \widebar{\bm{\psi}_{sr}}   \right).\label{Lin1_4}
\end{align}
By replacing \eqref{Lin1_3} and \eqref{Lin1_4} in objective function \eqref{B2S_StoPD_upperOF}, we can reformulate the expected profit of a stochastic power producer as the linear expression:
\begin{align}
& -\sum_{sr} \pi_{s} \widetilde{\bm{w}}^T_{sr}diag(\bm{q}^T\bm{U}) \widebar{\bm{\psi}_{sr}} - \bm{i}^T\bm{U}\bm{1} .\label{Lin1_5}
\end{align}
This way, we turn the original investment problem into a one-level mixed-integer linear programming problem that can be solved using commercial optimization software. 

To end this section, we compare the linear expressions \eqref{Lin2_5} and \eqref{Lin1_5} of the profit of the stochastic power producer under ConvMC and StocMC, respectively. Note that \eqref{Lin1_5} consists of two terms, one depending on the dual variable $\widebar{\bm{\psi}_{sr}}$, and the other representing the total investment cost. For given investment decisions $\bm{U}$, $\widebar{\bm{\psi}_{sr}}$ in \eqref{Bilevel2StageSto} indicates the sensitivity of the total system cost \eqref{B2S_Sto_lowerOF} with respect to a marginal increase of the realized stochastic power production $\widetilde{\bm{w}}^T_{sr}$, i.e., this dual variable captures the market value of each MWh of stochastic generation at each node of the system and scenario. Note that, if the stochastic power production cannot be completely injected in the system and spillage occurs, constraint \eqref{B2S_Sto_lowerC6} is non-binding and the value of its shadow price, $\widebar{\bm{\psi}_{sr}}$, is equal to 0. Therefore, expression \eqref{Lin1_5} implies that the investor aims at locating its stochastic generating units at the sites such that the difference between the expected market value of the stochastic resources \Big($-\sum_{sr} \pi_{s} \widetilde{\bm{w}}^T_{sr}diag(\bm{q}^T\bm{U}) \widebar{\bm{\psi}_{sr}}$\Big) and the investment cost of the undertaken projects \big($\bm{i}^T\bm{U}\bm{1}$\big) is maximized.

Conversely, the expression \eqref{Lin2_5} of the stochastic producer profit under ConvMC has an extra term, if compared with \eqref{Lin1_5}, namely, $\sum_s\pi_s\widehat{\bm{w}}^T_{s}diag(\bm{q}^T\bm{U})\big(\bm{\lambda}_{s}-\sum_r \bm{\lambda}^R_{sr} \big)$. This term arises from the fact that, under ConvMC, stochastic generating units are dispatched to their expected production in the day-ahead market, which is equivalent to say that stochastic producers take a short forward position of $\widehat{\bm{w}}_s$ MWh in this market. Consequently, under ConvMC, the power producer is not solely interested in locations with the highest expected market value of the stochastic energy resources, but also in those locations where the expected pay-off of this forward contract, computed as $\sum_s\pi_s\widehat{\bm{w}}^T_{s} diag(\bm{q}^T\bm{U})\big(\bm{\lambda}_{s}-\sum_r \bm{\lambda}^R_{sr} \big)$, is maximized. 

\section{Uncertainty characterization} \label{SectionUncertainty}

This section describes the procedure to generate the scenario set that characterizes the uncertainty pertaining to stochastic power production and demand both in the day-ahead and the balancing stages. The correlation between wind power production and load is assumed to be zero as it happens in the Nordic area \citep{holttinen2005impact}.

The uncertain parameters affecting the clearing of the day-ahead market are the conditional expectation of both the electricity consumption and stochastic production. Since no inter-temporal constraints, such as ramping limits or minimum times, are considered in the proposed market-clearing models, day-ahead scenarios $s$ are generated by randomly sampling the marginal probability distribution of each uncertain parameter. With respect to the demand, for the sake of simplicity, only one uncertain parameter characterizing the total system load is considered. The probability distribution of the hourly total load can be determined based on historical values \citep{murphy2005generation}. Accordingly, the load level corresponding to each node is computed as a fix percentage of the total consumption, denoted by $\chi_n$. 

In the same way, the probability distributions describing the stochastic power production at different locations can be obtained using historic data. For simplicity, the random variables corresponding to the stochastic power production at different locations are assumed to be independent and therefore, each day-ahead scenario is generated by randomly sampling each of the probability distributions separately. For example, if two wind farms are located at nodes $n_1$ and $n_2$, each day-ahead scenario $s$ consists of three values, i.e.,$\left(\widehat{L}^T_{s}, \widehat{W}_{n_1s}, \widehat{W}_{n_2s} \right)$. In order to capture the uncertain behaviour of all stochastic parameters, a large enough number of day-ahead scenarios needs to be generated. On the other hand, a scenario reduction technique can be applied in order to trim down the number of scenarios while keeping most of the stochastic information embedded in such scenarios \citep{dupavcova2003scenario, morales2009scenario}.

Once the day-ahead scenario set is defined, we also need to model the forecast errors of stochastic production and demand. These errors translate into energy deviations with respect to the day-ahead schedule that are to be settled in the balancing market. Therefore, for each day-ahead scenario $s$, we need to generate an additional set of scenarios $r$ that characterize the forecast errors of demand and stochastic power production. According to \cite{billinton1996reliability}, the load forecast error can be described by a normal distribution with a zero mean and a given standard deviation. Particularly, if we consider wind power as the only type of stochastic production in the system, we can model the wind forecast error using beta distributions with zero mean and a standard deviation that depends on the normalized predicted wind power production as proposed in \cite{Fabbri}. More specifically, they use empirical data to prove that the standard deviation of the wind power forecast error 24 hours in advance can be determined according to the following expression:
\begin{equation}
\sigma = \kappa(0.01837 + 0.20355 \cdot p), 
\label{StandDeviation}
\end{equation}
where $p$ represents the conditional expectation of the wind power in per unit, and $\kappa$ is a parameter that allows us to adjust the variability of the wind forecast error between perfect forecast ($\kappa = 0$) and the base case ($\kappa = 1$). Without loss of generality, we consider here that the day-ahead market for energy delivery throughout day $d$ is cleared at 12 p.m. the previous day. Consequently, the time delay between the market clearing and the actual operation of the system varies from 12 to 36 hours. In this paper, however, we have assumed that the probability distribution of the wind forecast error is determined in all cases for a time lag of 24 hours, which is the average look-ahead time.

Using the value of the standard deviation provided by \eqref{StandDeviation}, the shape parameters of the beta distribution $\alpha$ and $\beta$ modeling the wind power forecast error are computed as follows:
\begin{equation}
\alpha = \frac{\left(1-p\right)\cdot p\cdot p}{\sigma^2}-p, \quad \quad \quad \quad \beta = \alpha \left( \frac{1-p}{p}   \right).
\end{equation} 

\begin{figure}[htbp]
	\centering		\includegraphics[scale=0.4]{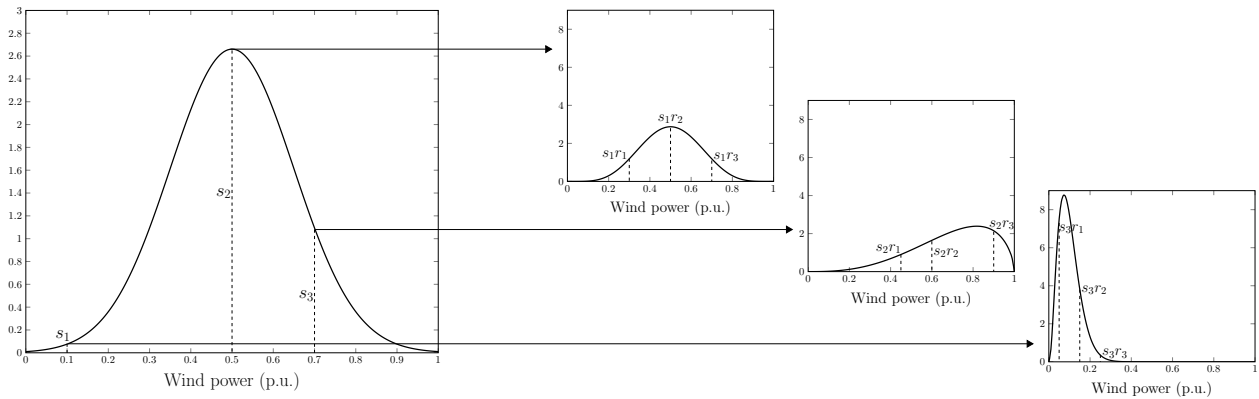}
	\caption{Uncertainty characterization of the day-ahead wind power production forecast (left-hand side figure) and the wind forecast error for particular day-ahead forecast values (three right-hand side figures)}
	\label{fig:UncertCharact}
\end{figure}

Once the probability distributions of the wind and demand forecast error are known for a particular day-ahead scenario $s$, a set of balancing scenarios $r$ is generated by randomly sampling these distributions as illustrated in Figure \ref{fig:UncertCharact}, where the left-hand side figure characterizes the hourly variation of the conditional expectation of the wind power production throughout the planning horizon, while each of the right-hand side plots represents the probability distribution of the wind forecast error linked to particular values of the wind forecast at the day-ahead stage. We also assume here that random variables characterizing forecasting errors of demand and wind power production are independent. Likewise, scenario reduction techniques can also be used to reduce the size of the scenario tree at the balancing stage. Following the previous example, each scenario $r$ at the this stage consists of three values: the total load forecasting error $\widetilde{L}^T_{sr}$, as well as the wind power production forecast error at $n_1$ and $n_2$, i.e., $\widetilde{W}_{n_1sr}$ and $\widetilde{W}_{n_2sr}$.

Apart from the uncertainty due to short-term variations of demand and stochastic power production, investment decisions also face other type of long-term uncertainties such as annual demand growth, fuel prices, network or generation expansion by other players, or regulatory risk. In this paper, however, we focus our attention on the impact of the short-term market design on investment decisions of stochastic producers and therefore, we solve next the proposed investment models for a given target year, in which the effect of long-term uncertainties can be disregarded. Nevertheless, it should be noted that investment models \eqref{Bilevel2StageConv} and \eqref{Bilevel2StageSto} can be expanded to the multi-year case, albeit increasing the computational burden of such models.

\section{Illustrative example} \label{SectionExample}

\begin{figure}[htbp]
	\centering		\includegraphics{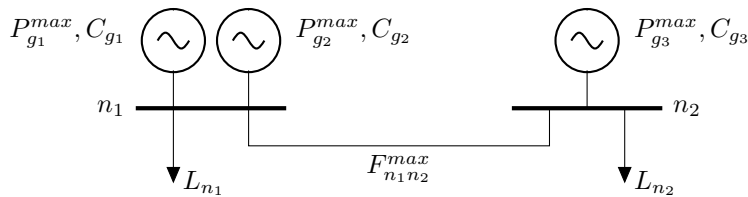}
	\caption{Two-node illustrative example}
	\label{fig:2node}
\end{figure}

The two-node system shown in Figure \ref{fig:2node} is used to illustrate the main features of the proposed investment models \eqref{Bilevel2StageConv} and \eqref{Bilevel2StageSto}. The data of the three thermal generating units are provided in Table \ref{table:2-node unit data}. Observe that $g_1$ and $g_3$ are, in comparative terms, cheaper than $g_2$, but inflexible. On the other hand,  the most expensive unit $g_2$ can increase or decrease its production up to 50MW at the balancing stage. Note that the price offers to sale or repurchase energy by unit $g_2$ in the balancing market slightly deviate from its price offer in the day-ahead market, indicating that producer $g_2$ is willing to be flexible in return for a price premium on the energy traded for balancing \citep{pritchard2010single}.

\begin{table}[htb] \centering 
\caption{Generating unit data of the two-node example. $P^{max}_{g}$ and $C_{g}$ represent the capacity and the marginal cost of unit $g$. Likewise, the pairs $(P^{max,u}_{g},C^{u}_{g})$, and $(P^{max,d}_{g},C^{d}_{g})$ are the upper bounds and prices for up and down balancing power provided by unit $g$, in that order. Powers and prices measured in MW and \$/MWh, respectively}
\begin{tabular}{c c c c c c c }    
  \hline
  \up\down{\vphantom{\Large{$A_p$}}} $g$ & $P^{max}_{g}$ & $C_{g}$ & $P^{max,u}_{g}$ & $C^{u}_{g}$ & $P^{max,d}_{g}$ & $C^{d}_{g}$ \\
  \hline
  \up$g_1$ & 400 & 20   & -  & -   & -  & -   \\  
  $g_2$ & 400 & 30   & 50 & 35  & 50 & 29  \\   
  \down$g_3$ & 600 & 22   & -  & -   & -  & -   \\
  \hline
\end{tabular}\label{table:2-node unit data}\end{table}

\begin{figure}[htbp]
\centering
\includegraphics[scale =0.5]{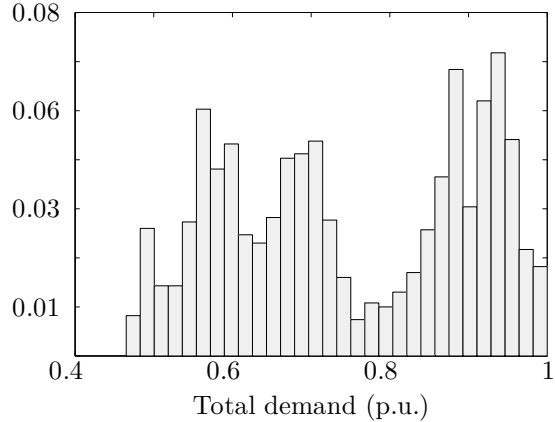}
\caption{Probability mass function of the per unit total demand}\label{fig:TotalLoad}
\end{figure}

Nodes $n_1$ and $n_2$ are connected through a 200-MW line with a susceptance of 68.5 p.u. and a base power of 100 MW. The variability of the day-ahead forecast for the total system demand is obtained by scaling the yearly load profile provided in \citep{grigg1999ieee}, and depicted in Figure \ref{fig:TotalLoad} in per unit, to a net peak demand of 660 MW, with the consumption level at $n_2$ being ten times bigger than the demand at $n_1$ in all cases. For illustration purposes, in this example, the errors associated with the day-ahead point forecast of demand are considered to be zero.

\begin{figure}[htbp]
\centering
\subfigure[Node $n_1$]{\includegraphics[scale =0.5] {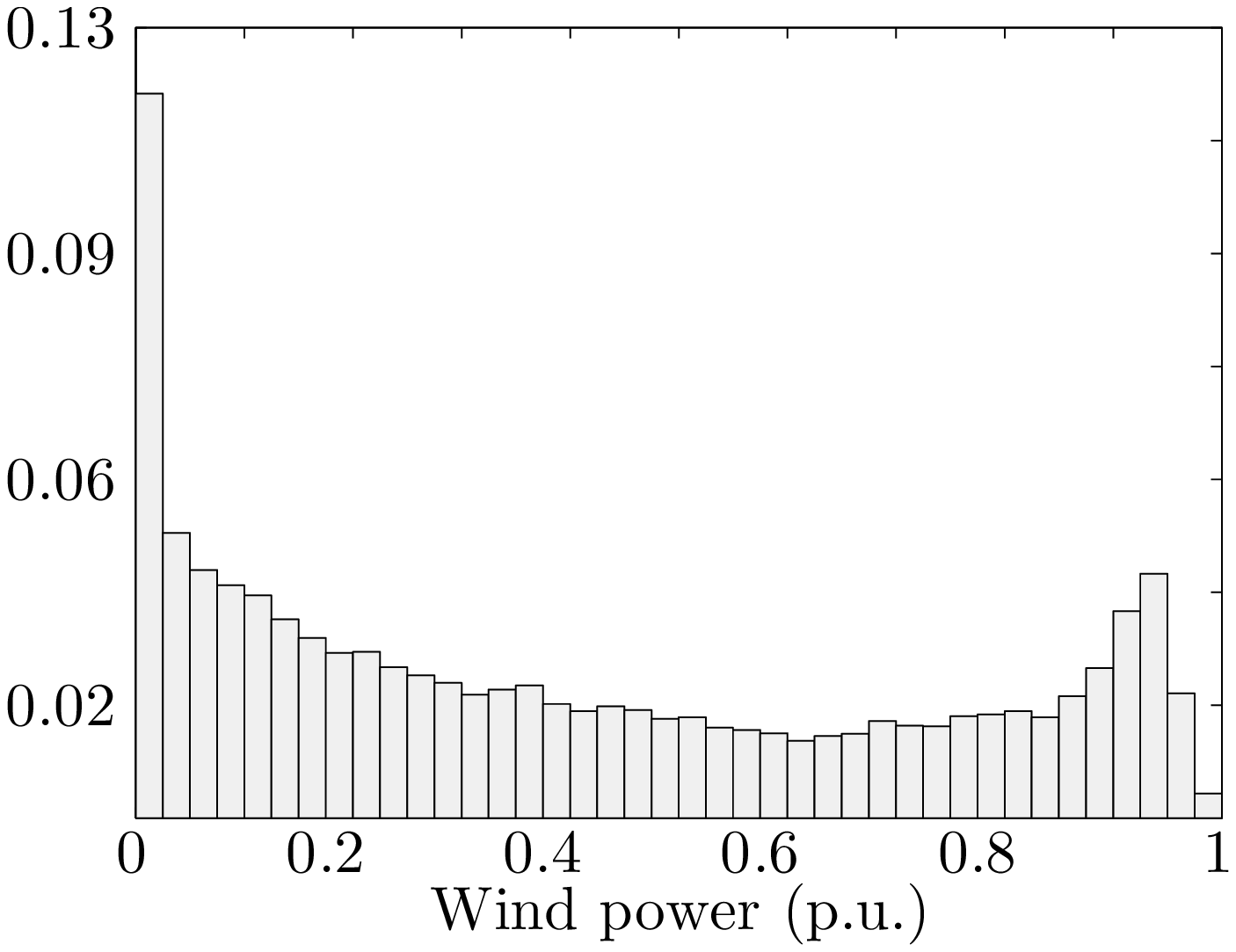}\label{fig:Windnode1}}
\subfigure[Node $n_2$]{\includegraphics[scale =0.5] {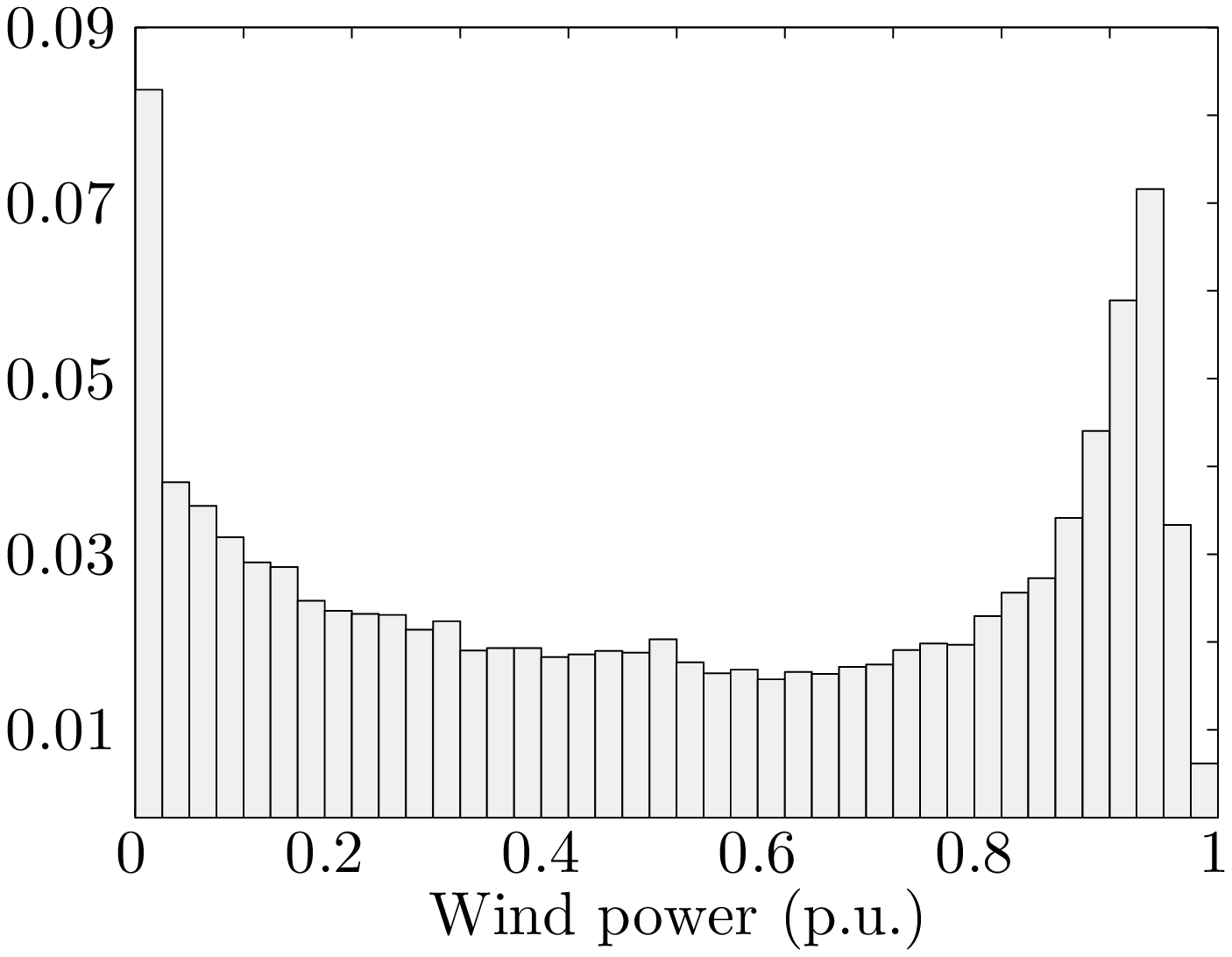}\label{fig:Windnode2}}
\caption{Probability mass function of the per unit wind power production}\label{fig:Windnode}
\end{figure}

Figures \ref{fig:Windnode1} and \ref{fig:Windnode2} show the frequency distribution for the power produced (normalized by the capacity of each wind farm) during the year 2006 by two simulated wind farms located at South Dakota (45$^{\circ}$13' North, 96$^{\circ}$55' West) and Montana (45$^{\circ}$20' North, 104$^{\circ}$24' West) in the US, respectively. These two distributions are used in this example to characterize the wind power production variability at the day-ahead stage at nodes $n_1$ and $n_2$. These data have been simulated by the National Renewable Energy Laboratory (NREL) of the US Department of Energy in the context of the Eastern Wind Integration and Transmission Study (EWITS). These data can be freely downloaded from \cite{NREL2010}. It is worth mentioning that the average wind power production at node $n_1$ (0.3960 p.u.) is lower than that at $n_2$ (0.4899 p.u.). Finally, the scenarios representing the forecast error distribution of each wind power production level are generated as described in Section \ref{SectionUncertainty}. For simplicity, the forecast errors at $n_1$ and $n_2$ are also assumed to be uncorrelated.

By randomly sampling the probability distributions depicted in Figures \ref{fig:TotalLoad} and \ref{fig:Windnode}, we generate 2000 scenarios characterizing the different situations of the day-ahead market throughout the planning horizon, each scenario being comprised of three values, namely, the per-unit total demand in the system and the wind power production at nodes $n_1$ and $n_2$. This initial set is then reduced to 100 scenarios using the technique described in \cite{dupavcova2003scenario}. For each of these day-ahead scenarios, we similarly obtain 1000 samples from the probability distributions of the wind forecast production error at nodes $n_1$ and $n_2$, and reduce them to 20 scenarios. Note that, since the demand forecast errors are considered zero in this example, each balancing scenario only contains two values, namely, the wind power production errors at nodes $n_1$ and $n_2$. The final scenario set is then composed of $100\times20$ scenarios.

In this example, we present two different analysis. In the first one, the producer has to decide the optimal location of a single 200-MW wind farm within the two-node system of Figure \ref{fig:2node}. In the second case, ten possible investment projects consisting of a 20-MW wind farm each are considered.

To investigate how each market-clearing mechanism, namely, ConvMC and StocMC, affects the location of the 200-MW wind farm, the main results of solving optimization problems \eqref{Bilevel2StageConvPD2} and \eqref{Bilevel2StageStoPD} are collated in Table \ref{table:Results example}. Wind power producer revenues and system costs are computed for a target year and expressed in million\$. Moreover, it is assumed that the power producer only owns the new wind farm to be installed. The investment cost of the wind farm is supposed to be proportional to its capacity with a rate of \$800/kW, and independent of its location. For an average useful life of 40 years, the annualized investment cost is determined as $\dfrac{200\cdot10^3\times800}{40} = \$4\cdot10^6$.

\begin{table}[htb]\begin{center}
\caption{Optimal investment decisions corresponding to the 200-MW wind farm for the two market clearing mechanisms (ConvMC and StocMC)}\label{table:Results example}
\begin{tabular}{l c c c}    
  \hline
  \up\down                 					& ConvMC  	&& StocMC  	\\
  \hline
  \up Optimal location 						& $n_1$  	&& $n_2$ 	\\  
   Day-ahead producer revenue (m\$)			& 14.0 		&& 12.7		\\
   Balancing producer revenue (m\$)			& -9.0		&& -2.4 	\\
   Annualized investment cost (m\$) 		& 4.0		&&	4.0		\\
   Total producer profit (m\$)     			& 1.0 		&& 6.3 		\\      
   Day-ahead system cost (m\$)  			& 77.2 		&& 81.7 	\\
   Balancing system cost (m\$)  			& 5.0 		&& -4.4		\\
   Total system cost (m\$)					& 82.2		&& 77.3 	\\
   \down Demand covered by wind (\%)		& 15.3		&& 18.2		\\ 
  \hline
\end{tabular}\end{center}\end{table}

First notice the significant difference of the balancing costs incurred by the power producer depending on the market-clearing design and investment decisions. While the most efficient use of the system flexibility made by StocMC results in relatively lower balancing costs for the wind power producer ($2.4$m\$), the day-ahead dispatch obtained by ConvMC entails balancing costs that are of the same order of magnitude as the revenue earned by the power producer at the day-ahead market ($9.0$m\$). Consequently, the profit made by the wind power producer is more than six times higher when investment decisions are made under a market design in which a prognosis of the future balancing costs is used to determine day-ahead dispatch decisions.

\begin{figure}[htbp]
\centering
\subfigure[Day-ahead wind dispatch]{\includegraphics[scale =0.55] {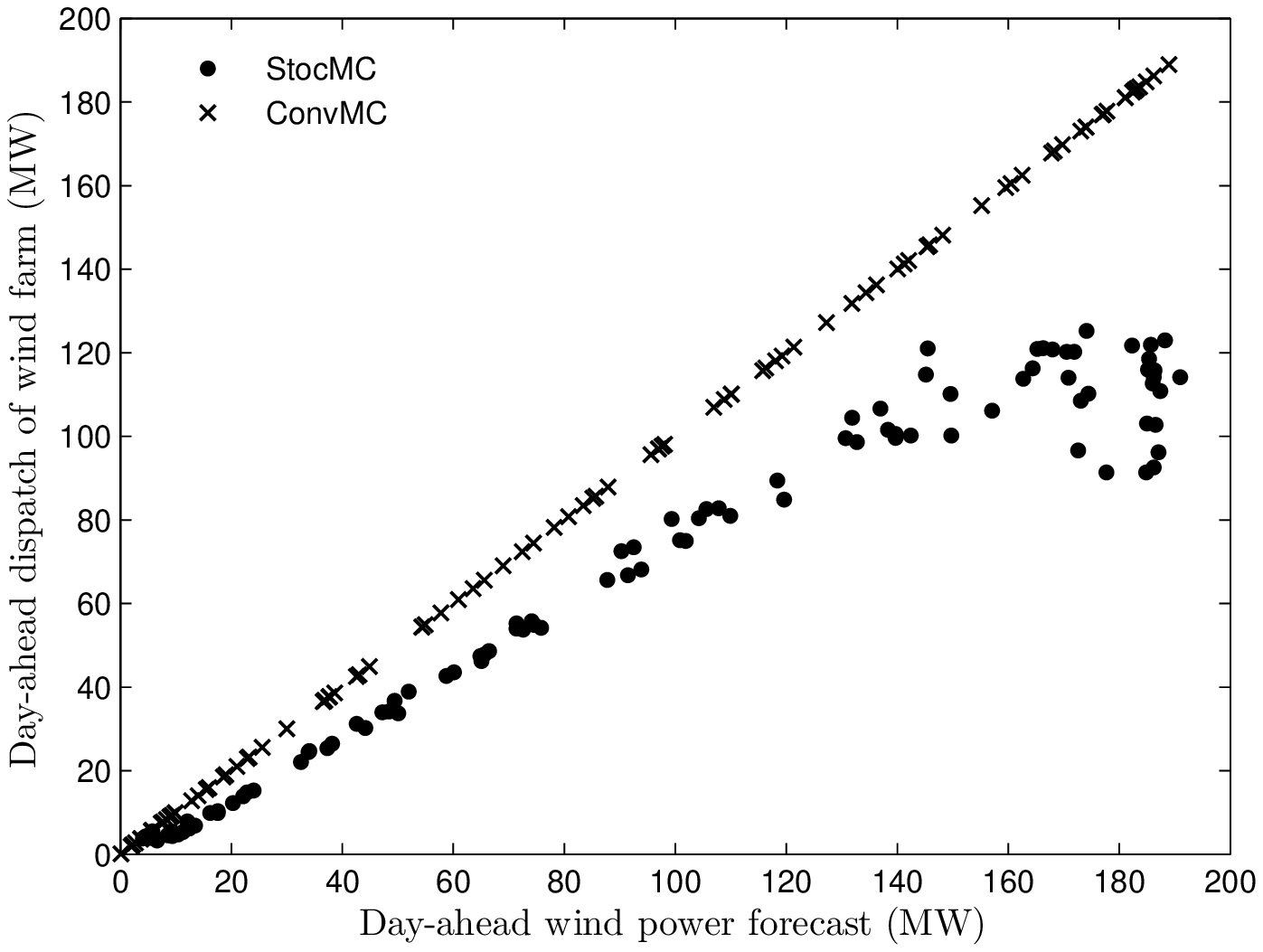}\label{fig:ExampleWindDispatch}}
\subfigure[Day-ahead dispatch of $g_2$]{\includegraphics[scale =0.55] {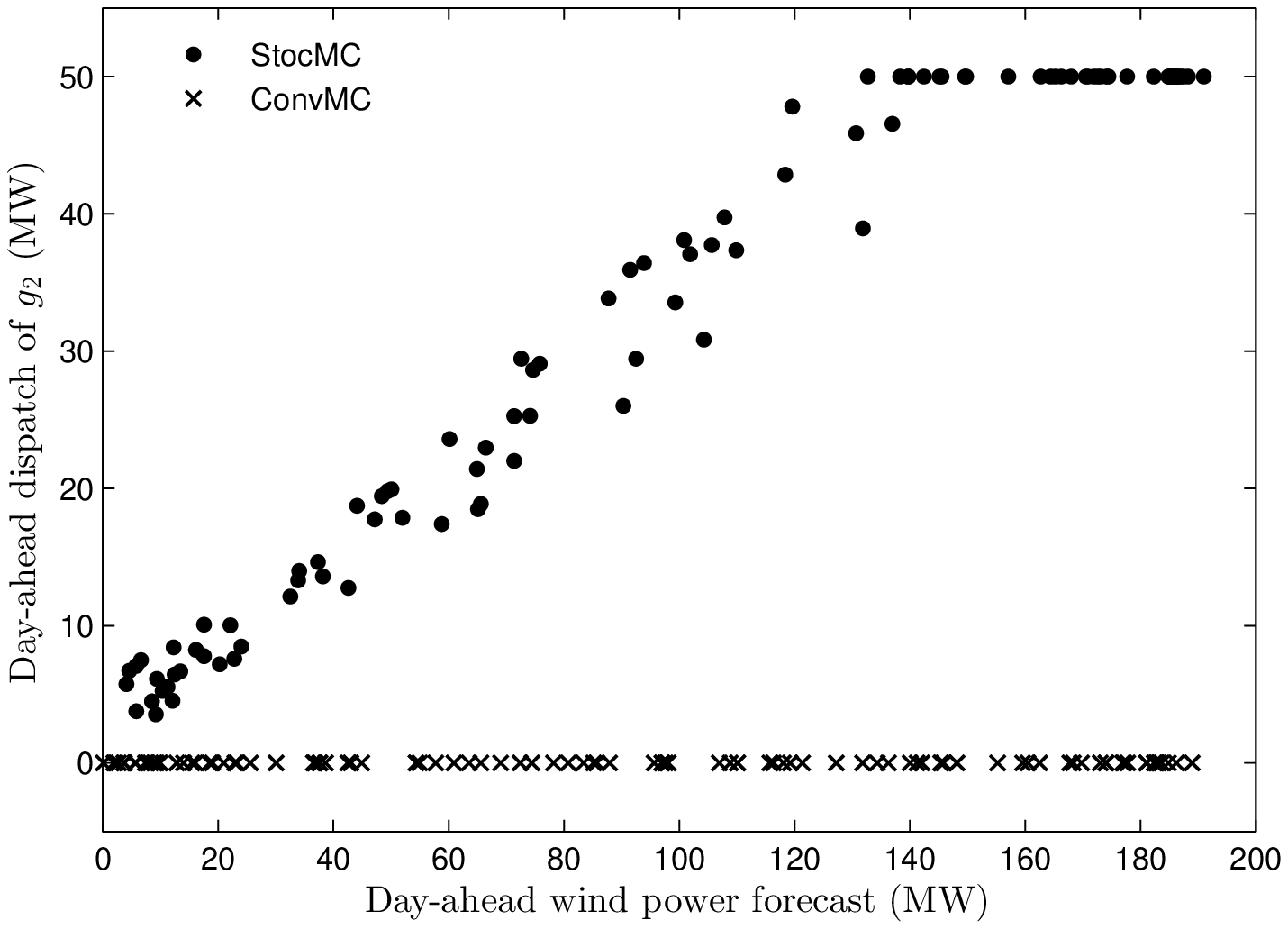}\label{fig:ExampleDispatchG2}}
\caption{Day-ahead dispatch of the 200-MW wind farm (a) and the flexible unit $g_2$ (b) as a function of the day-ahead wind power forecast and the two market-clearing mechanisms (ConvMC and StocMC)} \label{fig:Example1}
\end{figure}

In order to further illustrate the impact of the market-clearing design on wind investment decisions, we plot in Figure \ref{fig:Example1} the day-ahead dispatch of both the new wind farm and the flexible generating unit $g_2$ as a function of the wind power production forecast at the day-ahead stage for both market models. It can be observed that, for all wind power forecast levels, ConvMC dispatches the wind farm according to its forecast, and the generating unit $g_2$ to 0 MW for being the most expensive unit. By doing so, this market design attains the lowest day-ahead production cost (as indicated in Table \ref{table:Results example}), but fails to efficiently balance the system due to the lack of up- and down-regulating resources in those scenarios in which the actual wind production differs from the forecast value. On the other hand, StocMC dispatches the wind farm to a value lower than its forecast production, in particular, to a value such that the 50 MW of upward regulating capacity provided by unit $g_2$ are enough to mostly cover the lack of wind power generating at the balancing stage. Furthermore, StocMC dispatches unit $g_2$ out of cost-merit order to take full advantage of its relatively cheap downward regulation capacity in those scenarios where there is a surplus of wind power production. 

The particular features of these two paradigmatic market models have, in turn, a significant impact on the investment decisions in new wind generating capacity, as also shown in the results of Table \ref{table:Results example}. Indeed, even though lower wind resources are available at node $n_1$, the inefficient use of system flexibility under ConvMC leads the power producer to locate the 200-MW wind farm at this location in exchange for easy access to the balancing resources provided by unit $g_2$, also located at this node. On the other hand, the power producer places the 200-MW wind farm at node $n_2$ under StocMC to exploit the more favorable wind resources at this location, without significantly increasing its imbalance cost thanks to a better utilization of the flexible generating unit $g_2$ under this market design. Apart from the profit increase of the power producer, the impact of StocMC on the wind farm location proves also to be beneficial for the power system efficiency and the social welfare by increasing the percentage of demand covered by wind power production from 15.3\% to 18.2\% and reducing the total system cost by 6\%. In this sense, it should be noticed that the stochastic market-clearing procedure results in a lower operation cost not only because it makes a more efficient use the flexibility provided by $g_2$, but also because it facilitates the installation of the new wind farm at the node with the highest capacity factor by ensuring cheaper access to balancing resources at this location. 

To conclude this section we solve optimization models \eqref{Bilevel2StageConvPD2} and \eqref{Bilevel2StageStoPD} considering ten  20-MW wind farm projects in order to determine how many of them should be undertaken as well as their optimal locations. We assume that all the wind farms are composed of the same type of wind turbine, being the yearly investment cost of each of the projects equal to \$0.4m. The two-node system of Figure \ref{fig:2node} and the data formerly presented are also used in this analysis.

\begin{figure}[htbp]
\centering
\subfigure[ConvMC]{\includegraphics[scale =0.5] {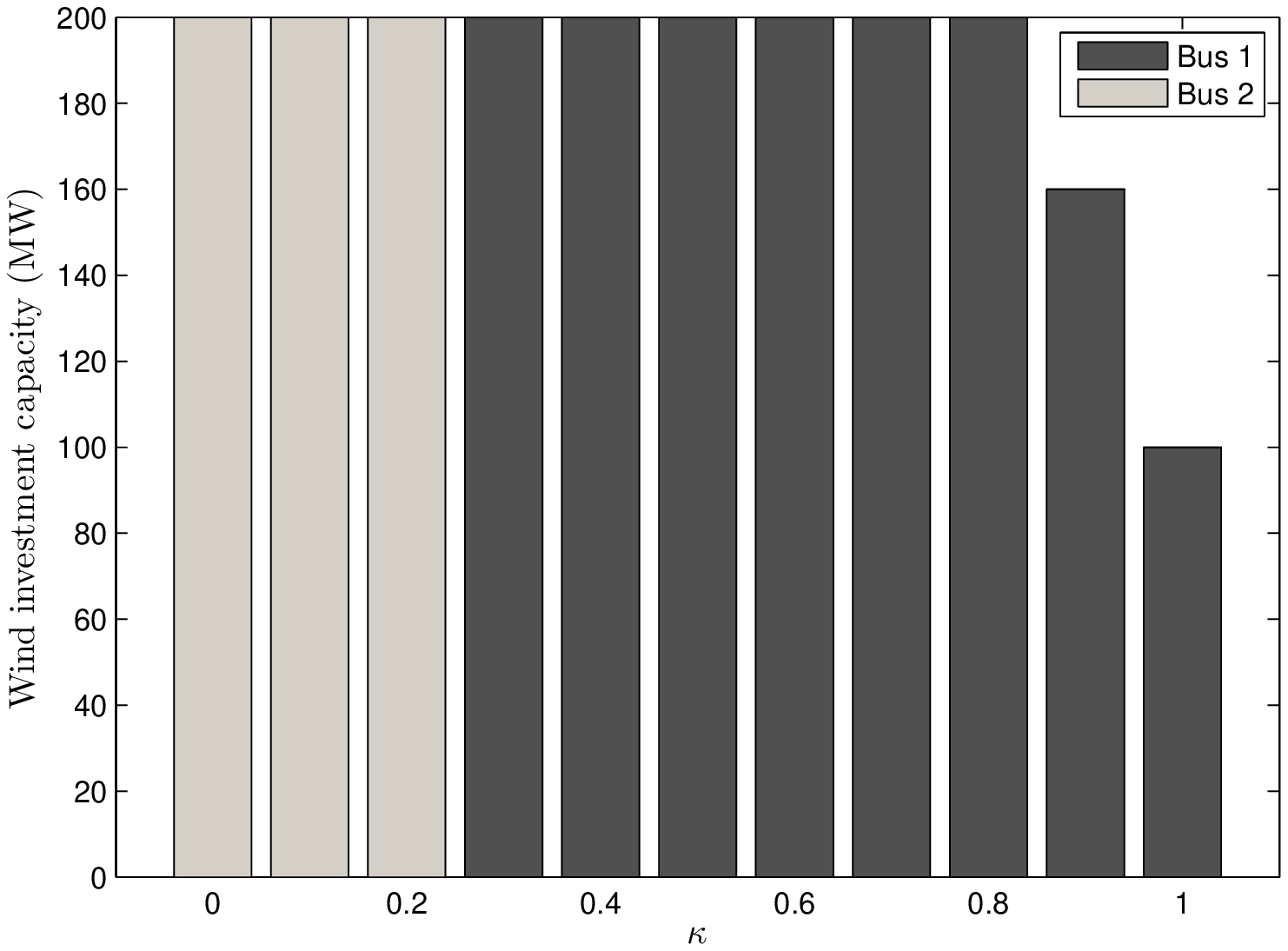}\label{fig:InvetExampleConvMC}}
\subfigure[StocMC]{\includegraphics[scale =0.5] {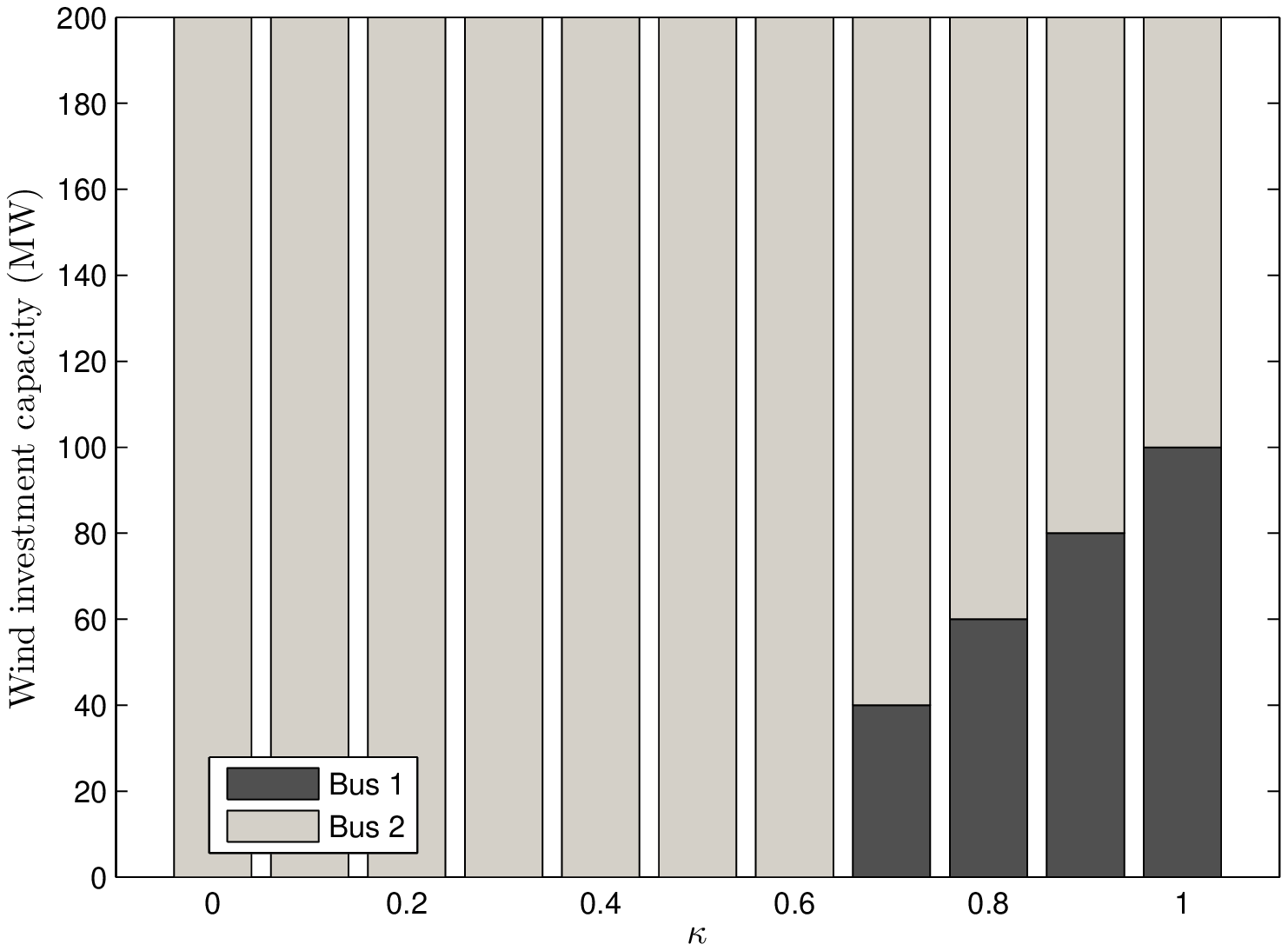}\label{fig:InvetExampleStocMC}}
\caption{Optimal investment decisions corresponding to the ten 20-MW wind farm projects as a function of the wind power forecast error ($\kappa$) and the two market-clearing mechanisms (ConvMC and StocMC)}\label{fig:InvestExample}
\end{figure}

Figure \ref{fig:InvestExample} depicts the optimal investment decisions of the power producer as a function of the parameter $\kappa$, which scales the standard deviation of the wind power production forecast error, for both the conventional and stochastic market-clearing mechanisms. Observe that, for $\kappa=0$, which implies that the actual wind power production is equal to the forecast value in all cases, the optimal investment decisions consists in installing the ten 20-MW wind farms at node $n_2$ in order to exploit the superior wind resources at this location. Since load forecast errors are disregarded, no balancing actions are required, and therefore, ConvMC and StocMC lead to the same investment decisions. 

As the value of $\kappa$ increases, the optimal investment decisions under each market-clearing design significantly diverge. In the case of ConvMC, Figure \ref{fig:InvetExampleConvMC} shows that for values of $\kappa$ higher than 0.2, it is more profitable to allocate all the wind farms at $n_1$, thus giving up on higher wind resources to get access to the cheaper balancing power provided by generating unit $g_2$, which is located at this node. Note also that, for the highest values of $\kappa$, this market-clearing design results in load shedding and wind spillage events, which dramatically increases the balancing cost of the wind power producer. In order to keep imbalance costs sufficiently low and ensure profitability, the power producer only invests in 100 MW of wind capacity for $\kappa = 1$. Conversely, the power producer installs all the wind farms at $n_2$ under StocMC, as long as $\kappa$ remains below 0.6 thanks to the ability of this market-clearing mechanism to deal with real-time imbalances in a more efficient manner. Although the power producer starts to increasingly relocate wind farms at node $n_1$ for values of $\kappa$ higher than 0.6, the optimal use of system flexibility achieved by StocMC allows the power producer to maintain the investment level at 200 MW even for $\kappa = 1$.

\begin{figure}[htbp]
\centering
\subfigure[Wind producer profit]{\includegraphics[scale =0.5] {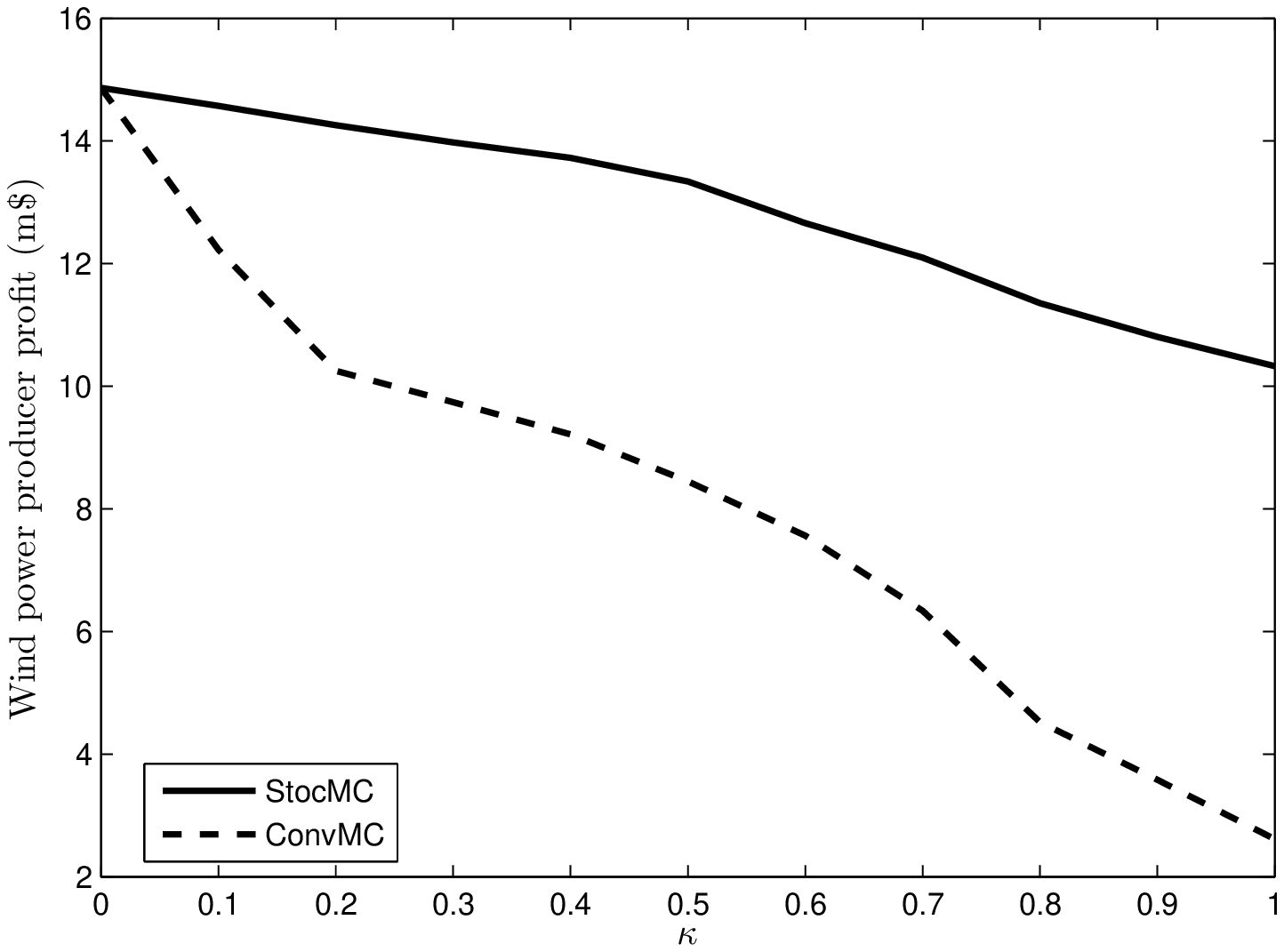}\label{fig:ProducerProfitExample}}
\subfigure[Wind share]{\includegraphics[scale =0.5] {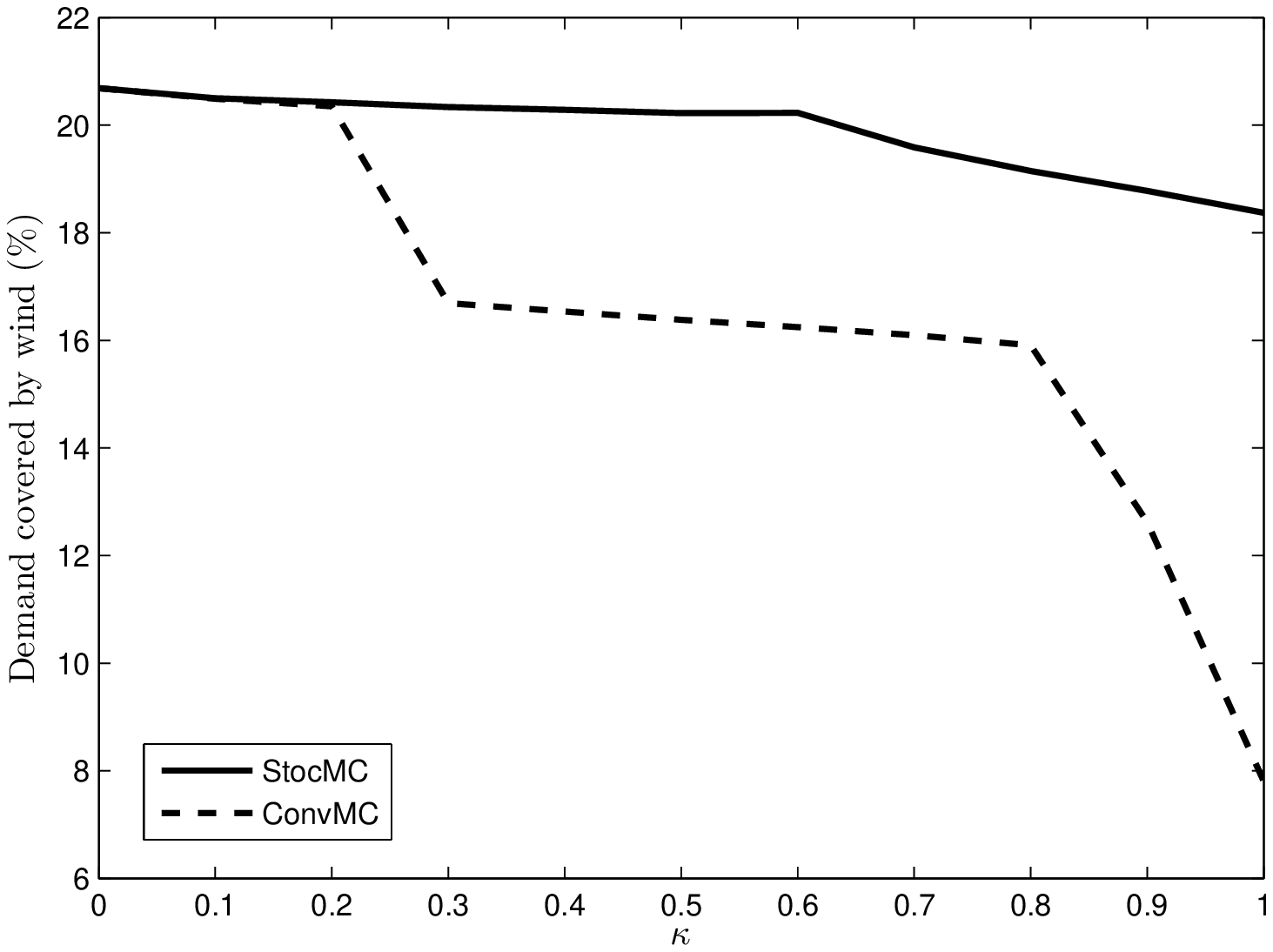}\label{fig:DemCovWindExample}}
\caption{Wind power producer profit and wind share as a function of the wind power forecast error ($\kappa$) and the two market-clearing mechanisms (ConvMC and StocMC) }\label{fig:ResultsInvestExample}
\end{figure}

For completeness, Figure \ref{fig:ResultsInvestExample} plots the evolution of the wind power producer profit and the percentage of demand covered by wind generation as a function of the parameter $\kappa$ for both market-clearing designs. Based on the depicted results, the following conclusions are in order. First, higher values of forecast errors involve a decrease of both the power producer profit and the share of wind generation in the electricity supply. Second, investment decisions under StocMC entail a significantly higher profit for investors in wind capacity, if compared to the profit obtained under ConvMC. For example, for $\kappa = 1$, the optimal investment decisions under ConvMC and StocMC involve a yearly profit of 2.61m\$ and 10.32m\$, respectively. Finally, investment decisions made under the stochastic market-clearing mechanism lead to higher percentages of demand covered by wind power production not only because it provides a more efficient use of the available balancing resources, but also because StocMC leads both to higher levels of installed wind capacity and to a better utilization of the available wind resources. Indeed, for $\kappa = 1$, this percentage is equal to 7.8\% and 18.4\% for ConvMC and StocMC, respectively. 

\section{24-node case study} \label{SectionCaseStudy}

\begin{figure}[htbp]
	\centering		\includegraphics[scale=0.6]{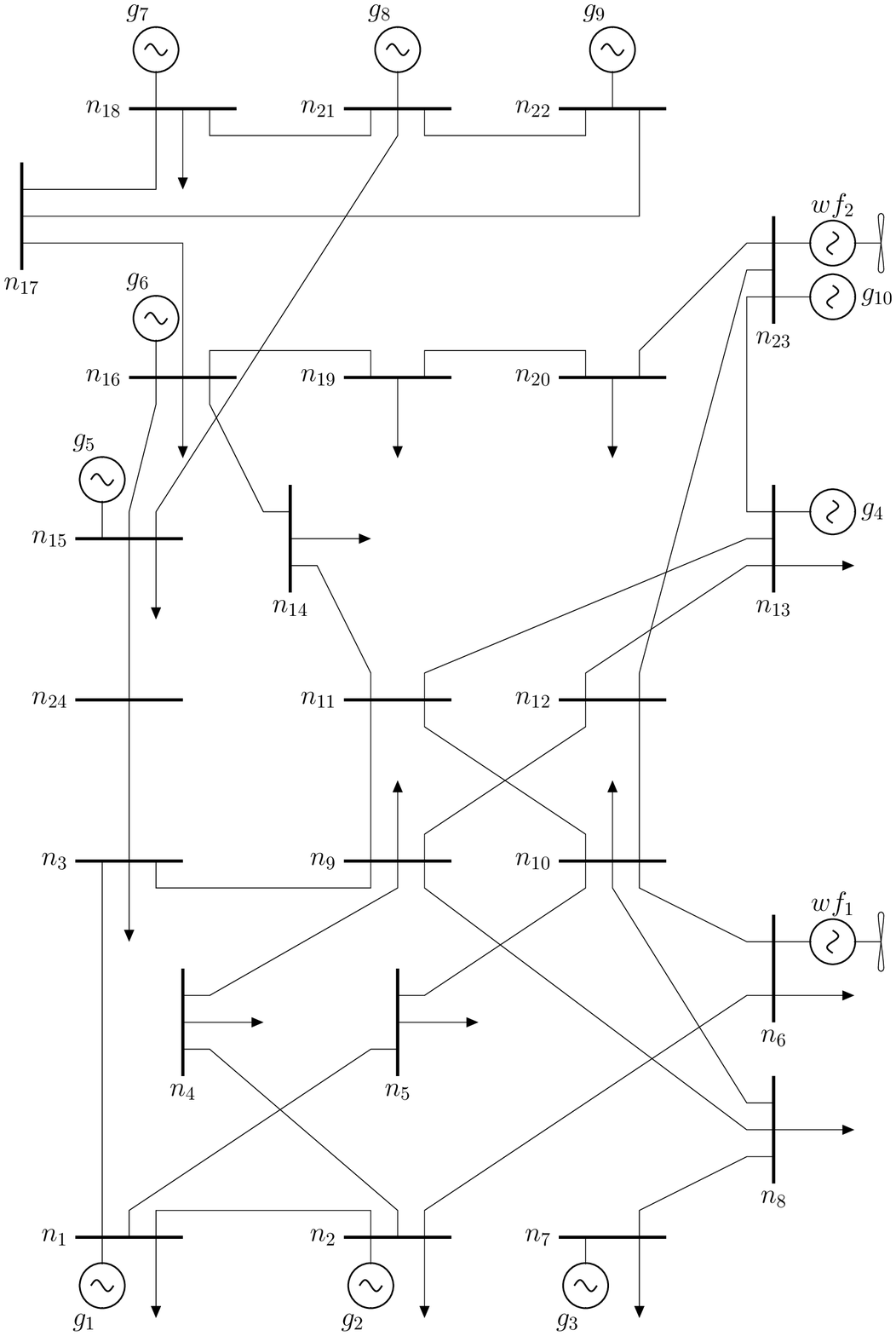}
	\caption{24-node system}
	\label{fig:24node}
\end{figure}

In this section the investment models \eqref{Bilevel2StageConvPD2} and \eqref{Bilevel2StageStoPD} are further discussed using a modified version of the 24-node system presented in \cite{grigg1999ieee} and depicted in Figure \ref{fig:24node}. The technical characteristics and offers for the day-ahead and the balancing market of the ten thermal generating units considered for this study are provided in Table \ref{table:24-node unit data}. Note that the most expensive generating units $g_7$, $g_8$, and $g_9$ are flexible and fast, and offer their entire power capacities for balancing at relatively competitive prices. On the other hand, the rest of the existing thermal generating units are assumed to be inflexible and only offer $\pm10$ MW of up/down balancing power at $\pm100$\$/MWh as a representation of the automatic generation control reserves of these units. Furthermore, two 100-MW wind farms are located at nodes $n_6$ and $n_{23}$. Finally, Table \ref{table:linedata} provides the characteristics of the network.

\begin{table}[htb]\begin{center}
\caption{Generating unit data of the 24-node case study. $P^{max}_{g}$ and $C_{g}$ represent the capacity and the marginal cost of unit $g$. Likewise, the pairs $(P^{max,u}_{g},C^{u}_{g})$, and $(P^{max,d}_{g},C^{d}_{g})$ are the upper bounds and prices for up and down balancing power provided by unit $g$. Powers and prices measured in MW and \$/MWh, respectively}\label{table:24-node unit data}
\begin{tabular}{c c c c c c c c c c c c c c c c c }    
  \hline
  \up\down{\vphantom{\Large{$A_p$}}} $g$ & $n$ & $P^{max}_{g}$ & $C_{g}$ & $P^{max,u}_{g}$ & $C^{u}_{g}$ & $P^{max,d}_{g}$ & $C^{d}_{g}$ && $g$ & $n$ & $P^{max}_{g}$ & $C_{g}$ & $P^{max,u}_{g}$ & $C^{u}_{g}$ & $P^{max,d}_{g}$ & $C^{d}_{g}$\\
  \hline
\up $g_1$    & $n_1$    & 400 & 25.9   & 10    & 100     & 10    & -100   && $g_6$    & $n_{16}$ & 400 & 19.2   & 10    & 100     & 10    & -100   \\  
$g_2$    & $n_2$    & 575 & 22.3   & 10    & 100     & 10    & -100   && 
$g_7$    & $n_{18}$ & 120 & 30.1   & 120   & 35.1    & 120   & 29.1   \\  
$g_3$    & $n_7$    & 500 & 26.6   & 10    & 100     & 10    & -100   &&  
$g_8$    & $n_{21}$ & 100 & 30.6   & 100   & 35.6    & 100   & 29.6   \\  
$g_4$    & $n_{13}$ & 520 & 21.2   & 10    & 100     & 10    & -100   &&  
$g_9$    & $n_{22}$ & 80  & 31.1   & 80    & 36.1    & 80    & 30.1   \\  
$g_5$    & $n_{15}$ & 475 & 17.5   & 10    & 100     & 10    & -100   && 
\down $g_{10}$ & $n_{23}$ & 450 & 20.8   & 10    & 100     & 10    & -100   \\    
    \hline
\end{tabular}\end{center}\end{table}

\begin{table}[htb]\begin{center}
\caption{Network data of the 24-node case study. $B_{nm}$ and $F^{max}_{nm}$ represent the susceptance and the capacity of the transmission line connecting nodes $n$ and $m$, measured in p.u. and MW, respectively}\label{table:linedata}
\begin{tabular}{c c c c c c c c c c c c c c c}
  \hline
  \up\down $nm$ 	& $B_{nm}$  & $F^{max}_{nm}$  & & $nm$ 	& $B_{nm}$  & $F^{max}_{nm}$  & & $nm$   & $B_{nm}$  & $F^{max}_{nm}$ & & $nm$   & $B_{nm}$  & $F^{max}_{nm}$ \\
  \hline  
  \up $n_{1}n_{2}$	  & 68.5 & 800  & &	$n_{6}n_{10}$	  & 15.6  & 400  & & $n_{11}n_{14}$  & 23.5	& 500 & & $n_{16}n_{17}$ & 38.0	& 500\\
  $n_{1}n_{3}$	  & 4.4  & 400  & & 	$n_{7}n_{8}$	  & 15.3  & 400  & & $n_{12}n_{13}$  & 20.5	& 500 & & $n_{16}n_{19}$ & 42.7	& 500\\
  $n_{1}n_{5}$	  & 11.0 & 400& &	$n_{8}n_{9}$	  & 5.7   & 400  & &$n_{13}n_{23}$  & 10.2	& 500 & & $n_{17}n_{18}$ & 69.9	& 500\\
  $n_{2}n_{4}$	  & 7.4  & 400  & &	$n_{8}n_{10}$   & 5.7	  & 400  & &$n_{13}n_{23}$  & 11.3	& 500 & & $n_{17}n_{22}$ & 9.4	  & 500\\
  $n_{2}n_{6}$	  & 4.9  & 400 &  &	$n_{9}n_{11}$   & 11.9	& 100  & &$n_{14}n_{16}$  & 16.8	& 500 &  &$n_{18}n_{21}$ & 75.8	& 1000\\
  $n_{3}n_{9}$	  & 7.9  & 400 &  &	$n_{9}n_{12}$   & 11.9	& 100  & &$n_{15}n_{16}$  & 58.1	& 500 & & $n_{19}n_{20}$ & 49.3	& 1000\\
  $n_{3}n_{24}$	  & 11.9 & 100 &  &	$n_{10}n_{11}$  & 11.9	& 100  & &$n_{15}n_{21}$  & 40.2	& 1000 & & $n_{20}n_{23}$ & 89.3	& 1000\\
  $n_{4}n_{9}$	  & 9.0  & 400  &  &	$n_{10}n_{12}$  & 11.9	& 100  & &$n_{15}n_{24}$  & 18.9	& 500 & & $n_{21}n_{22}$ & 14.5	& 500\\
  \down $n_{5}n_{10}$	  & 10.6 & 400  &  &	$n_{11}n_{13}$  & 20.5	& 500  &  &  & 	& 	& &  \\ 
  \hline
\end{tabular}\end{center}\end{table}

We determine now, for the two market-clearing mechanisms, the optimal generating expansion decisions pertaining to four 250-MW wind farms that can only be placed at nodes $n_1$, $n_2$, $n_{18}$, and $n_{21}$. The investment cost is assumed to be proportional to the capacity of each wind farm at a rate of \$800/kW and independent of its location, being the payback period equal to 40 years.

For simplicity, the conditional expectation of the wind power production at nodes $n_1$, $n_2$ and $n_6$ is assumed to be perfectly correlated and characterized by the probability mass function plotted in Figure \ref{fig:Windnode1}. Similarly, Figure \ref{fig:Windnode2} contains the expected wind power production at nodes $n_{18}$, $n_{21}$ and $n_{23}$. Furthermore, the wind power production at the two group of nodes is considered to be uncorrelated. Lastly, the probability distribution of the wind power forecast error at each node is determined as explained in Section \ref{SectionUncertainty}. 

The hourly total system load is described in \cite{grigg1999ieee} and plotted in Figure \ref{fig:TotalLoad} in p.u. in the form of a probability mass function. In this case study, a peak load of 2850 MW is considered. Moreover, Table \ref{table:PercentageLoad} contains the percentages of the system load corresponding to the consumption at each node. The forecast error of the total load is modeled as a normal distribution with zero mean and a standard deviation of 2\% of the mean value. Likewise, the demand forecast error is proportionally distributed among all nodes according to the parameter $\chi_n$. For simplicity, the demand is assumed to be independent of the wind power production.

\begin{table}[htb]\begin{center}
\caption{Percentage of system load of the 24-node case study}\label{table:PercentageLoad}
\begin{tabular}{c c c c c c c c c c c c c c c c c c}
  \hline
  \up $n$ & $n_{1}$ & $n_{2}$& $n_{3}$& $n_{4}$& $n_{5}$& $n_{6}$& $n_{7}$ & $n_{8}$ & $n_{9}$ & $n_{10}$ & $n_{13}$& $n_{14}$& $n_{15}$& $n_{16}$& $n_{18}$& $n_{19}$& $n_{20}$ \\
  \down $\chi_n (\%)$ & 3.8 & 3.4 & 6.3 & 2.6 & 2.5 & 4.8 & 4.4 & 6.0 & 6.1 & 6.8 & 9.3 & 6.8 & 11.1 & 3.5 & 11.7 & 6.4 & 4.5 \\  
  \hline
\end{tabular}\end{center}\end{table}

The variability of the conditional expectation of the wind power and the total system demand at the day-ahead stage is characterized by a scenario set of 2000 scenarios, which is reduced to a final set of 30 scenarios. For each of these day-ahead scenarios, 1000 scenarios are generated and reduced to 20 to model the forecast errors of both parameters at the balancing stage. As in the example presented in Section \ref{SectionExample}, the variability of the forecast error is artificially modified via the parameter $\kappa$ for the sake of analysis.

Table \ref{table:InvCaseStudy} shows the optimal investment decisions for different values of $\kappa$ and the two market-clearing models corresponding to optimization problems \eqref{Bilevel2StageConvPD2} and \eqref{Bilevel2StageStoPD}. We can observe that if the wind power production of all wind farms in the system can be perfectly forecast ($\kappa=0$), the investor decides to locate the four wind farms at node $n_2$ due to the higher availability of wind resources at this node. Since no wind forecast errors need to be settled at the balancing stage, both ConvMC and StocMC provide the same optimal investment decisions and the same profit for the wind power producer.

\begin{table}[htb]\begin{center}
\caption{Investment decisions in the 24-node case study}\label{table:InvCaseStudy}
\begin{tabular}{l c c c c c c c}
\hline
\up\down &  \multicolumn{3}{c}{ConvMC}   && \multicolumn{3}{c}{StocMC}   \\
\hline
\up &  $\kappa=0$ & $\kappa=0.5$ &$\kappa=1$ &&$\kappa=0$ &$\kappa=0.5$ &$\kappa=1$ \\
Locations  & $n_2,n_2,n_2,n_2$  &$n_2,n_{18},n_{18},n_{21}$  &$n_{18},n_{21}$   && $n_2,n_2,n_2,n_2 $  & $n_1,n_1,n_1,n_{18} $ & $n_1,n_1,n_{18},n_{18}$ \\  
Profit  (m\$)     &  72.3   & 41.8   & 13.2   && 72.3   & 63.4   & 42.2  \\  
Total cost (m\$) &  281.5  & 313.8  & 354.4  && 280.2  & 291.8  & 309.5  \\  
\down Wind penet. (\%)  &  27.4   & 23.0   & 12.4   && 27.5   & 26.0   & 23.4  \\  
  \hline
\end{tabular}\end{center}\end{table}

Conversely, if the wind forecast errors are different from zero ($\kappa=0.5$), it can be observed that, under ConvMC, the maximum profit is attained if three out of the four wind farms are placed at nodes $n_{18}$ and $n_{21}$, where cheaper balancing resources are available at the expense of decreasing the capacity factor of these units. Note that these decisions do not only have negative consequences for the power producer, whose profit is almost halved, but also for the whole system by decreasing the share of wind production in the electricity supply and increasing consequently the total production cost. On the other hand, the StocMC and its more efficient day-ahead dispatch of wind and thermal generating units allows the wind power producer to maintain three wind farms at node $n_1$ to take advantage of  higher capacity factors, while still limiting the imbalance costs due to wind forecast errors. 

Finally, let us consider the case in which $\kappa=1$. Under ConvMC, the imbalance costs become so high that the best option for the power producer is to build only two of the four available projects, thus dramatically reducing the wind share in the system to 12.4\%. On the other hand, under StocMC, it is still profitable to invest in the four available projects, two of which are to be located at node $n_{18}$ with lower wind resources but easier access to the cheaper balancing resources provided by units $g_7$, $g_8$, and $g_9$. 

Results in Table \ref{table:InvCaseStudy} also show that the consideration of wind forecast errors when making wind investment decisions in this 24-node system reduces the profit of the wind producer by 81.7\% and 41.6\% depending on whether day-ahead dispatch decisions are made according to ConvMC or StocMC, respectively. Besides, it can be observed that the wind power producer is willing to invest in higher amounts of wind power capacity under  StocMC, which results in a wind share in the electricity supply twice as high as that under ConvMC.

\section{Conclusions} \label{SectionConclusions}

The impossibility of accurately forecasting the power production of stochastic generating units at the clearing of the day-ahead market implies that a significant share of the profit of stochastic power producers comes from the trading of their production imbalances in the balancing market. In this paper, we present a mathematical framework to model the impact of these imbalance costs on generation expansion decisions pertaining to stochastic generating units. The proposed models are formulated as MPEC that determine the investment decisions maximizing the profit of the stochastic producer including the revenues from both the day-ahead and balancing markets. Day-ahead and balancing prices are endogenously generated through a set of lower-level problems that represent the different market conditions throughout the planning horizon. Uncertainty pertaining to stochastic production and demand level is modeled via scenarios. The impact of two paradigmatic market-clearing mechanisms on investment decisions is also analyzed and discussed. The investment models proposed in this paper are recast as mixed-integer linear programming problems that can be solved using commercial optimization software. 

Two main conclusions can be drawn from the results presented throughout this paper. First, it is shown that imbalance costs incurred by stochastic power producers may significantly affect optimal generating investment decisions. In this sense, it can be more profitable to install new stochastic generating units at locations with a lower capacity factor, but cheaper access to balancing resources. Secondly, a market-clearing design that makes a more efficient use of system flexibility and provides, therefore, cheaper balancing power to stochastic producers will, in the long run, promote higher amounts of installed stochastic capacity in those locations with plentiful renewable resources. In turn, this market design will lead to higher penetration levels of stochastic electricity production into power systems and lower operating costs. 

\ACKNOWLEDGMENT{
S. Pineda thanks the Danish Council for Strategic Research for support
through “5s - Future Electricity Markets” project, no. 12-
132636/DSF.

J. M. Morales thanks the Danish Council for Strategic Research for support
through ENSYMORA project, no. 10-093904/DSF. J. M. Morales is also sponsored by DONG Energy

The authors would also like to thank M. Zugno and T. Boomsma for
insightful comments and relevant observations about the models and results presented in this paper.
}

\bibliographystyle{ormsv080}

\end{document}